\documentclass[oneside,12pt,dvips]{amsart}
\usepackage[dvipdf]{graphicx}
\usepackage{euler, amsfonts, amssymb, amsmath, latexsym, epsfig,epic,hyperref,pstricks}
\usepackage[dvips]{geometry}

\DeclareFontFamily{T1}{calligra}{}
\DeclareFontShape{T1}{calligra}{m}{n} {<-> callig15}{}

\newtheorem{thm}{Theorem}[section]
\newtheorem{prop}[thm]{Proposition}
\newtheorem{cor}[thm]{Corollary}
\newtheorem{lem}[thm]{Lemma}

\theoremstyle{definition}
\newtheorem*{defn}{Definition}
\newtheorem*{conj}{Conjecture}

\theoremstyle{remark}

\newtheorem{rem}{\it Remark}[section]
\newtheorem{ex}{\it Example}[section]

\numberwithin{equation}{section} 

\newcommand{\bT}{\mathbf{Tab}}
\newcommand{\bR}{\mathbf R}
\newcommand{\bS}{\mathbf{S}}
\newcommand{\bB}{\mathbf{B}}
\newcommand{\bG}{\mathbf{G}}
\newcommand{\pr}{^{\prime}}
\newcommand{\prpr}{^{\prime\prime}}
\newcommand{\codim}{{\rm codim\,}}
\newcommand{\Rank}{{\rm Rank\,}}
\newcommand{\sh}{{\rm sh\,}}
\newcommand{\ov}{\overline}

\newcommand{\B}{\mathcal{B}}
\newcommand{\Bb}{\mathfrak{B}}
\newcommand{\T}{\mathcal{T}}

\newcommand{\X}{\mathcal{X}}

\newcommand{\V}{\mathcal{V}}
\newcommand{\F}{\mathcal{F}}

\newcommand{\Or}{\mathcal{O}}

\newcommand{\G}{\mathfrak{g}}
\newcommand{\Sl}{\mathfrak{s}\mathfrak{l}}
\newcommand{\Gl}{\mathfrak{g}\mathfrak{l}}
\newcommand{\gh}{\mathfrak{h}}

\newcommand{\bo}{\mathfrak{b}}

\newcommand{\nil}{\mathfrak{n}}

\newcommand{\Ze}{\mathbb Z}
\newcommand{\uar}{\uparrow}
\newcommand{\dar}{\downarrow}
\newcommand{\rar}{\rightarrow}
\newcommand{\lar}{\leftarrow}
\begin{document}


\title[two column case]{Intersections of components of a Springer fiber
of codimension one for the two column case}

\author{A. Melnikov}
\address{Department of Mathematics, University of Haifa, Haifa
31905, Israel.} \email{melnikov@math.haifa.ac.il}

\author{N.G.J. Pagnon}
\address{Chennai Mathematical Institute
Plot H1, SIPCOT IT Park Padur PO, Siruseri 603103 India.}
\email{pagnon@cmi.ac.in}


\begin{abstract}
This paper is a subsequent paper of \cite{MPI}. Here we consider the
irreducible components of Springer fibres (or orbital varieties) for
two-column case in ${\rm GL}_n$. We describe the intersection of two
irreducible components, and specially give the necessary and
sufficient condition for this intersection to be of codimension one.
This work find some motivations in the conjecture 6.3 of Kazhdan and
Lusztig \cite{KL} or else in the work of Hotta \cite{Hot}. Since an
orbital variety in two-column case is a finite union of the Borel
orbits, we solve the initial question by determining orbits of
codimension one in the closure of a given orbit. We show that they
are parameterized by a specific set of involutions called {\em
descendants}, already introduced by the first author in a previous
work. Applying this result we show that the intersections of two
components of codimension one are irreducible and provide the
combinatorial description in terms of Young tableaux of the pairs of
such components.
\end{abstract}
\maketitle

\section{Introduction}
\subsection{}\label{1.1}
This paper is a continuation of \cite{MPI} and we adopt its
notation.

Let $\bG=GL_n$ and respectively $\G=\Gl_n={\rm Lie}(\bG)$ on which
$\bG$ acts by conjugation. For $g\in \bG$ and $u\in\G$ we denote
this action by $g.u:=gug^{-1}.$

We fix the standard triangular decomposition
$\G=\nil\oplus\gh\oplus\nil^-$ where $\nil$ is the subalgebra of
strictly upper triangular matrices, $\nil^-$ is the subalgebra of
strictly lower triangular matrices and $\gh$ is the subalgebra of
diagonal matrices. The associated Weyl group in this case is
identified with the symmetric group $\bS_n.$ Let
$\bo:=\gh\oplus\nil$ be the standard Borel subalgebra and  $\bB$ the
Borel subgroup of $\bG$ with ${\rm Lie}(\bB)=\bo$ that is the
subgroup of invertible upper-triangular matrices.
 For $x\in {\nil}$ let $\Or_x=\bG.x$ be its orbit. Consider $\Or_x\cap{\nil}.$
Its irreducible components are called {\bf orbital varieties}
associated to $\Or_x.$

Let $\F:=\bG/\bB$ be the flag manifold. For $x\in\nil$ put
$${\F}_{x}:=f^{-1}(x)=\{g\bB\in \F\ |\ x\in g.{\nil}\ \}= \{g\bB\in \F\
|\ g^{-1}.x\in {\nil}\ \} $$
The variety ${\F}_{x}$ is called the {\bf Springer fiber} above $x$.

 By Spaltenstein's construction \cite{Spa} there is a one to
one correspondence between irreducible components of $\F_x$ and
orbital varieties associated to $\Or_x.$ For $x\in\nil$ its Jordan
form is completely defined by a partition $\lambda=(\lambda_1,\ldots
\lambda_k)$ of $n$ where $\lambda_i$ is the length of $i-$th Jordan
block. Arrange the numbers of a partition $\lambda=(\lambda_1,\ldots
\lambda_k)$ in the decreasing order (that is
$\lambda_1\geq\lambda_2\geq\cdots\geq \lambda_k\geq 1$) and write
$J(x)=\lambda.$ Obviously $\Or_x$ and $\F_x$ are completely defined
by $J(x).$

In turn an ordered partition can be presented as a Young diagram
$D_{\lambda}$ -- an array with $k$ rows of boxes starting on the
left with the $i$-th row containing $\lambda_i$ boxes. In such a way
there is a bijection between Springer fibers (resp. nilpotent
orbits) and Young diagrams.

Fill the boxes of Young diagram $D_\lambda$ with $n$ distinct
positive integers. \ If the entries increase in rows from left to
right and in columns from top to bottom we call such an array a
Young tableau or simply a tableau of shape $\lambda.$ Let
$\bT_\lambda$ denote the set of all tableaux of shape $\lambda.$

Given $x\in\nil$ such that $J(x)=\lambda$ by Spaltenstein
(\cite{Spa0}) and Steinberg (\cite{Ste2}) there is a bijection
between components of $\F_x$ (resp. orbital varieties associated to
$\Or_\lambda$) and $\bT_\lambda$ (cf. \ref{2.3}). For $T\in
\bT_{\lambda}$ set $\F_T$ to be the corresponding component of
$\F_x.$ Respectively set $\V_T$ to be the corresponding orbital
variety associated to $\Or_\lambda.$

 Moreover by \cite{M-P} the number of
irreducible components and their codimensions in $\F_T\cap\F_{T\pr}$
is equal to the number of irreducible components and their
codimensions in $\V_T\cap \V_{T\pr}.$ We consider here the
intersections of orbital varieties.

\subsection{}\label{1.2}
In this paper we consider  orbital varieties of nilpotent order 2.
They correspond to tableaux with two columns so this case is called
two column case.

 For the convenience we use the conjugate partitions. For
$x\in\nil$ of nilpotent order 2 that is such that $J(x)=(2,\ldots)$
we put $\sh(x):=(n-k,k)$ if the corresponding $J(x)$ consists of $k$
blocks of length 2 and $n-2k$ blocks of length 1. In other words
$(n-k,k)$ is the conjugate partition of $(2,\ldots).$ Respectively
we put $\sh(T):=(n-k,k)$ if the corresponding Young diagram has the
first column of  length $n-k$ and the second column of length $k.$
Respectively we put $\bT_{(n-k,k)}:=\{T\ :\ sh(T)=(n-k,k)\}.$

 We show that if $\codim_{\V_T}(\V_T\cap\V_S)=1$
for some $T,S\in \bT_{(n-k,k)}$ then $\V_T\cap\V_S$ is irreducible.
Further in Sections \ref{5.2a}, \ref{5.4} we describe
combinatorially such pairs of tableaux.
\subsection{}\label{1.3} Let us explain the results in some detail.
In our constructions we use intensively \cite{Mel1} and we adopt
partially its notation.
 Set $\X^2:=\{x\in\nil\ :\
x^{\scriptscriptstyle 2}=0\}$ to be the variety of nilpotent
upper-triangular $n\times n$ matrices of nilpotent order 2. For
$x\in \X^2$ set $\B_x:=\bB.x$ to be its $\bB$-orbit. Obviously
Jordan form is the same for all elements of $\B_x$ so we put
$\sh(\B_x):=\sh(x).$ Let $\Bb^2$ be the variety of $\bB$-orbits in
$\X^2.$ It is stratified by ranks. We put $\Bb_{(n-k,k)}:=\{\B\in
\Bb^2\ :\ \sh(\B)=(n-k,k)\}$ to be the subset of $\Bb^2$ of
$\bB$-orbits of rank $k.$

Set $\bS_n^2:=\{\sigma\in \bS_n\ :\ \sigma^2=Id\}.$   Let
$\bS_n^2(k)$ be the subset of the involutions containing $k$
disjoint 2-cycles. As it was shown in \cite{Mx2} there is a natural
bijection $\psi:\bS_n^2\rar \Bb^2$ such that for
$\sigma\in\bS_n^2(k)$ one has $\psi(\sigma)\in \Bb_{(n-k,k)}.$ Put
$\B_\sigma:=\psi(\sigma).$ As it is shown in \cite{Msmith} for any
$T\in\bT_{(n-k,k)}$ there is $\sigma_T\in \bS_n^2(k)$ such that
$\ov\B_{\sigma_T}=\ov\V_T.$ Moreover, by \cite{Mx2} these are all
the orbits of dimension $k(n-k)$ in $\Bb_{(n-k,k)}$ which is the
maximal possible dimension for an orbit in $\Bb_{(n-k,k)}.$

In \cite{Mel1} the closure of $\B_\sigma$ is described in terms of
involutions. The corresponding partial order on involutions is
determined by $\sigma'\preceq\sigma$ if
$\ov\B_{\sigma'}\subseteq\ov\B_\sigma$. In particular, given
$\sigma\in\bS_n^2(k)$ put

$$D(\sigma):=\{\sigma'\in\bS_n^2(k)\ :\ \sigma'\prec\sigma\ {\rm and}\
\sigma'\preceq \sigma\prpr\preceq \sigma\ \Rightarrow\
\sigma\prpr=\sigma\ {\rm or}\ \sigma\prpr=\sigma'\} $$
 This set (which we call the set of descendants of a given $\sigma$)
is described in \cite{Mel1}. Developing the results of \cite{Mel1}
we show in Theorem  \ref{thm3} that

$$D(\sigma)=\{\sigma'\in \bS_n^2(k)\ :\
\codim_{\ov\B_\sigma}\B_{\sigma'}=1\}$$
 This is the main technical
result of our paper.

As we show in \ref{4.9}  for $\sigma,\sigma'\in \bS_n^2(k)$ such
that $\dim\B_\sigma=\dim\B_{\sigma'}$ and
$\codim_{\ov\B_\sigma}(\ov\B_\sigma\cap\ov\B_{\sigma'})=1$ the
intersection $\ov\B_\sigma\cap\ov\B_{\sigma'}$ is reducible in
general. However as we show in \ref{5.1} for
$\sigma\prpr\in\bS_n^2(k)$ such that $\dim\B_{\sigma\prpr}=k(n-k)-1$
there exist exactly two involutions $\sigma,\sigma'\in\bS_n^2(k)$
such that $\sigma\prpr\in D(\sigma),D(\sigma')$ and moreover
$D(\sigma)\cap D(\sigma')=\{\sigma\prpr\}.$ As a corollary of this
result we get that the intersections of codimension 1 of orbital
varieties of nilpotent order 2 are irreducible.

For the sake of combinatorial completeness we study in \ref{4.9}
$\bS_n^2$ as a poset. In particular completing Theorem \ref{thm3} to
the cover of $\sigma$ (we call it the set of generalized
descendants) we show in Theorem \ref{lem2} that  $\sigma'$ is in the
cover of $\sigma$ iff $\codim_{\ov\B_\sigma}\ov\B_{\sigma'}=1$.

\subsection{}\label{1.4}
The body of the paper consists of 3 sections. In section 2 we give
the preliminaries in some detail and set the notation to make the
paper as self content as possible. In section 3 we study in detail
$\bB$-orbits of nilpotent order 2 and show that for
$\B_\sigma,\B_{\sigma'}\in\Bb_{(n-k,k)}$
$\codim_{\ov\B_\sigma}\B_{\sigma'}=1$ if and only if $\sigma'\in
D(\sigma).$ Finally in section 4 we apply the results of section 3
to orbital varieties and show that the intersections of orbital
varieties of codimension 1 are irreducible and give the description
of such orbital varieties in terms of Young tableaux.

At the end of the paper one can find the index of notation in which
symbols appearing frequently are given with the subsection where
they are defined. We hope that this will help the reader to find his
way through the paper.

\section{Preliminaries}

\subsection{}\label{2.1}
 For $\sigma\in\bS_n^2$ set $N_\sigma$ to be the ``strictly
upper-triangular part" of its permutation matrix, that is

\begin{eqnarray}\label{eq1}
(N_\sigma)_{i,j}:=\left\{\begin{array}{ll} 1&{\rm if}\ i<j\ {\rm
and}\
\sigma(i)=j;\\
0&{\rm otherwise.}\\
\end{array}\right.
\end{eqnarray}

\noindent Set  $\B_\sigma:=\B_{N_\sigma}$. For $\sigma\in\bS_n^2$
put $L(\sigma)$ to be the number of disjoint 2-cycles in it, let us
call the number $L(\sigma)$ the {\bf length} of the involution
$\sigma$ [do not confuse this notion with the usual definition given
for the minimal number of simple reflections in the writing of an
element of a Coxeter group]. In other words $L(\sigma)=k$ for any
$\sigma\in \bS_n^2(k).$ Note that $\sigma\in \bS_n^2(k)$ iff $\Rank
(\B_\sigma)=k$ that is iff $\B_\sigma\in\Bb_{(n-k,k)}.$ By
\cite[2.2]{Mx2} one has

\begin{prop}
$\Bb^2=\coprod\limits_{\sigma\in\bS_n^2}\B_{\sigma}.$
\end{prop}

\noindent In particular, for $k\ :\ 0\leq k\leq \frac{n}{2}$ one has
$\Bb_{(n-k,k)}=\coprod\limits_{\sigma\in\bS_n^2(k)}\B_{\sigma}.$

\subsection{}\label{2.2}
For $T\in \bT_{(n-k,k)}$ set $T=(T_1,T_2),$ where
$T_1=\left(\begin{array}{c}a_1\cr\vdots\cr
a_{n-k}\cr\end{array}\right)$ is the first column of $T$ and
$T_2=\left(\begin{array}{c}j_1\cr\vdots\cr j_k\cr\end{array}\right)$
is the second column of $T.$ It is enough to define the columns as
sets since the entries increase from  top  to  bottom in the
columns. We denote $\langle T_i\rangle$ when we consider it as a
set.

Put $\sigma_{\scriptscriptstyle T}=(i_1,j_1)\ldots (i_k,j_k)$ where
$i_1=j_1-1,$ and $i_s=\max\{d\in \langle
T_1\rangle\setminus\{i_1,\ldots,i_{s-1}\}\ |\ d<j_s\}$ for any
$s>1.$ For example, take
$$T=\begin{array}{ll}
1&4\\
2&5\\
3&7\\
6&8\\
\end{array}$$
Then $\sigma_{\scriptscriptstyle T}=(3,4)(2,5)(6,7)(1,8).$

Put $\B_T:=\B_{\sigma_T}.$ As it was shown in \cite[4.13]{Msmith}
\begin{prop}
For $T\in \bT_n^2$ one has $\overline\V_T=\overline\B_T.$
\end{prop}

\subsection{}\label{2.3}
 Write $\sigma\in \bS_n^2(k)$ as a product of disjoint cycles
of length 2. Order entries inside a given cycle in the increasing
order. Order the cycles in the increasing order according to the
first entry. Thus, $\sigma=(i_1,j_1)(i_2,j_2)\ldots(i_k,j_k)$ where
$i_s<j_s$ for any $1\leq s\leq k$ and $i_s<i_{s+1}$ for any $1\leq
s<k.$ We call this form the canonical form of $\sigma.$
\par

Given $\sigma=(i_1,j_1)(i_2,j_2)\ldots(i_k,j_k)$ where $i_s<j_s$ for
any $s\ :\ 1\leq s\leq k.$ For $s\ : 1\leq s\leq k$ set

$$q_s(\sigma):=q_{(i_s,j_s)}(\sigma):=\#\{i_p<i_s\ : \ j_p<j_s\}+
\#\{j_p\ :\ j_p<i_s\}. \eqno{(*)}
$$

Note that the definition of $q_s(\sigma)$ is independent of writing
$\sigma$ in the canonical form. However if it is written in the
canonical form then $q_1(\sigma)=0$ and to compute $q_s(\sigma)$ it
is enough to check only the pairs $(i_p,j_p)$ where $p<s.$

\begin{ex} Take $\sigma=(1,6)(3,4)(5,7).$ Then $L(\sigma)=3$ and
$q_1(\sigma)=0,\ q_2(\sigma)=0,\ q_3(\sigma)=2+1=3.$
\end{ex}
\par
By \cite[3.1]{Mx2} one has
\begin{thm}\label{dim}
 For $\sigma=(i_1,j_1)(i_2,j_2)\ldots(i_k,j_k)\in \bS_n^2(k)$ one has
$$\dim \B_\sigma=kn-\sum\limits_{s=1}^k(j_s-i_s)-\sum\limits_{s=2}^k q_s(\sigma).$$
\end{thm}
By \cite[Remark 5.5]{MPI} for $\sigma\in \bS_n^2(k)$ one has $\dim
\B_\sigma\leq k(n-k)$ and the equality is satisfied iff
$\sigma=\sigma_{\scriptscriptstyle T}$ for some $T\in\bT_{(n-k,k)}$.
Moreover if $\dim \B_\sigma= k(n-k)$ where
$\sigma=(i_1,j_1)\ldots(i_k,j_k)$ then $\sigma=\sigma_T$ where
$\langle T_2\rangle=\{j_1,\ldots, j_k\}$ and respectively $\langle
T_1\rangle=\{i\}_{i=1}^n\setminus\langle T_2\rangle.$
\subsection{}\label{2.4}
In \cite{Mel1} the combinatorial description of $\overline\B_\sigma$
(with respect to Zariski topology) for $\sigma\in\bS_n^2$ is
provided. Let us formulate this result.

For $1\leq i<j\leq n$ consider the canonical projections
$\pi_{i,j}:\nil_n\rightarrow \nil_{j-i+1}$ acting on a matrix by
deleting the first $i-1$ columns and rows and the last $n-j$ columns
and rows and define the rank matrix $R_u$ of $u\in\nil$ to be

\begin{eqnarray}\label{eq2}
(R_u)_{i,j}=\left\{ \begin{array}{ll}
0&{\rm if}\ i\geq j;\\
                       {\rm Rank}\,(\pi_{i,j}(u))&{\rm otherwise}.\\
                       \end{array}\right.
\end{eqnarray}

\noindent Obviously for any $y\in\B_u$ one has $R_y=R_u$ so that we
can define $R_{\B_u}:=R_u.$ Put $R_\sigma:=R_{\B_\sigma}.$

Let $\Ze^+$ be the set of non-negative integers. Put
$\bR^2_n:=\{R_\sigma\ :\ \sigma\in \bS_n^2\}.$ By \cite[3.1,
3.3]{Mel1} one has
\begin{prop} $R\in M_{n\times n}(\Ze^+)$  belongs to
$\bR^2_n$ if and only if it satisfies
\begin{itemize}
\item[(i)] $R_{i,j}=0$ if\ \ $i\geq j;$
\item[(ii)] For $i<j$ one has $R_{i+1,j}\leq R_{i,j}\leq R_{i+1,j}+1$
            and $R_{i,j-1}\leq R_{i,j}\leq R_{i,j-1}+1;$
\item[(iii)] If $R_{i,j}=R_{i+1,j}+1=R_{i,j-1}+1=R_{i+1,j-1}+1$ then
\begin{itemize}
\item[(a)] $R_{i,k}=R_{i+1,k}$ for any $k<j$ and
$R_{i,k}=R_{i+1,k}+1$ for any $k\geq j;$
\item[(b)] $R_{k,j}=R_{k,j-1}$ for any $k>i$ and
$R_{k,j}=R_{k,j-1}+1$ for any $k\leq i;$
\item[(c)] $R_{j,k}=R_{j+1,k}$ and $R_{k,i}=R_{k,i-1}$
for any $k\ :\ 1\leq k\leq n.$
\end{itemize}
\end{itemize}
\end{prop}
\subsection{}\label{2.5}
 Define the
following partial order on $M_{n\times n}(\Ze^+).$ For $A,B\in
M_{n\times n}(\Ze^+)$ put $A\preceq B$ if for any $i,j\ :\ 1\leq
i,j\leq n$ one has $A_{i,j}\leq B_{i,j}.$

The restriction of this order  to $\bR_n^2$ induces a partial order
on $\bS_n^2$ by setting $\sigma\succeq \sigma\pr$ if
$R_{\sigma}\succeq R_{\sigma'}$ for $\sigma,\sigma'\in \bS_n^2.$ By
\cite[3.5]{Mel1} one has

\begin{thm} For any $\sigma\in \bS_n^2$ one has
$$\overline\B_\sigma=\coprod\limits_{\sigma\pr\preceq\sigma}\B_{\sigma\pr}.$$
\end{thm}

\subsection{}\label{2.6}
Given $\sigma,\sigma\pr\in \bS_n^2$ we define $R_{\sigma,\sigma\pr}$
by
$$(R_{\sigma,\sigma\pr})_{i,j}:=\min\{(R_\sigma)_{i,j},(R_{\sigma\pr})_{i,j}\}.$$
Theorem 5.15 from \cite{MPI} claims

\begin{thm}
For any $\sigma,\sigma\pr\in \bS_n^2(k)$ one has
$$\overline
\B_\sigma\cap\overline\B_{\sigma\pr}=\coprod\limits_{R_{\sigma"}\preceq
R_{\sigma,\sigma'}}\B_{\sigma"}.$$
 In particular, $\overline\B_\sigma\cap\overline\B_{\sigma\pr}$ is
 irreducible if and only if $R_{\sigma,\sigma\pr}\in \bR_n^2.$ In
 that case
 $\overline\B_\sigma\cap\overline\B_{\sigma'}=\overline\B_\tau$
 where $R_\tau=R_{\sigma,\sigma\pr}.$
\end{thm}

\section{ $\bB$-orbits of nilpotent order 2}
\subsection{}\label{4.4}
In this section we make a more subtle analysis of the structure of
$\overline\B_\sigma.$

 Given
$\sigma=(i_1,j_1)\ldots(i_k,j_k)\in\bS_n^2(k)$ for $s:\ 1\leq s\leq
k$ set $\sigma_{(i_s,j_s)}^-$ to be the involution in $\bS_n^2(k-1)$
obtained from $\sigma$ by omitting the pair $(i_s,j_s).$

 In what
follows we need the following very simple lemma which is a
straightforward corollary of \cite[Lemma 5.10]{MPI}.

\begin{lem}
Let $\sigma=(i_1,j_1)\ldots(i_k,j_k).$ Then for any $s\ :\ 1\leq
s\leq k$ one has $\sigma=\sigma_{(i_s,j_s)}^-\cdot(i_s,j_s)$ and
$$(R_\sigma)_{i,j}=\left\{\begin{array}{lc}
(R_{\sigma_{(i_s,j_s)}^-})_{i,j} &{\rm if}\ i>i_s\ {\rm or}\ j<j_s;\\
(R_{\sigma_{(i_s,j_s)}^-})_{i,j}+1 & {\rm if}\ i\leq i_s\ {\rm and}\
j\geq j_s.\\
\end{array}\right.$$
In particular if $\sigma'=\sigma_{(i_s,j_s)}^-\cdot (i'_s,j'_s)$
where either $i'_s<i_s$ or $j'_s>j_s$ (or both) then
$\sigma'<\sigma.$
\end{lem}

\subsection{}\label{4.5}
Let $\sigma$ be some involution in $\bS_n(k).$ Set $D(\sigma)$ to be
the set of all $\sigma\pr\in \bS_n(k)$ such that

\begin{itemize}
\item[(a)] $\sigma\pr\prec\sigma;$
\item[(b)] for any $\sigma\prpr$  such
that $\sigma'\preceq \sigma\prpr\preceq \sigma$  one\ has either
$\sigma\prpr=\sigma$ or $\sigma\prpr=\sigma\pr.$
\end{itemize}

We call $D(\sigma)$ the set of {\bf descendants}. This set is
constructed in \cite[3.10-3.14]{Mel1} in order to prove Theorem
3.18. It contains of four types of elements. We give its description
in what follows.

To consider all the intersections of codimension 1 we also need to
define $A(\sigma)$ the set of {\bf ancestors} for
$\sigma\in\bS_n^2(k).$ Set $A(\sigma)$ to be the set of all
$\sigma'\in\bS_n^2(k)$ such that
\begin{itemize}
\item[(a)] $\sigma'\succ\sigma;$
\item[(c)] for any $\sigma\prpr$  such
that $\sigma'\succeq \sigma\prpr\succeq \sigma$  one\ has either
$\sigma\prpr=\sigma'$ or $\sigma\prpr=\sigma.$
\end{itemize}
We get this set from \cite[3.10-3.14]{Mel1}. It contains of four
types of elements exactly as $D(\sigma).$
\begin{rem}\label{rem2}
\begin{itemize}
\item[(i)]
Note that  for each $\lambda=(n-k,k)$ one has  $k\leq \frac{n}{2}$
so that $k<n-k+1$ and $\sigma_o(k):=(1,n-k+1)(2,n-k+2)\ldots(k,n)$
is defined. By \cite[Proposition 5.15]{MPI}
$\B_{(1,n-k+1)\ldots(k,n)}$ is the unique minimal orbit of
$\Bb_{(n-k,k)}.$ One has
$\dim\B_{(1,n-k+1)\ldots(k,n)}=\frac{1}{2}k(k+1)$ and
$\B_{(1,n-k+1)\ldots(k,n)}\subset\overline\B_{\sigma'}$ for any
$\sigma'$ such that $L(\sigma')\geq k.$ As well this is also the
only element of $\bS_n^2(k)$ such that its set of descendants is
empty.
\item[(ii)]
For any $\lambda=(n-k,k)$ there are exactly $\#\bT_{(n-k,k)}$ orbits
of  maximal dimension (equal $k(n-k)$) in $\Bb_{(n-k,k)}$ and these
correspond to the $B$-orbits $\B_T$, with $T\in\bT_{(n-k,k)}.$ As
well $\sigma_T$ (for $T\in \bT_{(n-k,k)}$) are the only elements of
$\bS_n^2(k)$ such that their set of ancestors is empty.
\end{itemize}
\end{rem}

\subsection{}\label{4.6a}
 To define the descendants and ancestors we will give another
description of the relation corresponding to the order $\preceq$ in
${\bf R}_n^{2}$.

 Given $\sigma=(i_1,j_1)\ldots(i_k,j_k)\in \bS_n^2$
set $\langle \sigma\rangle:=\{i_s,j_s\}_{s=1}^k$ to be its support
and $\langle \sigma\rangle^c:=\{s\}_{s=1}^n\setminus\langle
\sigma\rangle$ to be the compliment of its support.

For $i,j\ :\ 1\leq i<j\leq n$ put $\pi_{i,j}(\sigma)$ to be all the
pair of $\sigma$ with entries between $i$ and $j$, or formally
$\pi_{i,j}(\sigma)=(i_{s_1},j_{s_1})\ldots(i_{s_r},j_{s_r})$ where
for any $t\ :\ 1\leq t\leq r$ one has $i\leq i_{s_t},j_{s_t}\leq j.$
Note that this $\pi_{i,j}$ corresponds to the canonical projection
$\pi_{i,j}$ defined in \ref{2.4}, that is
$N_{\pi_{i,j}(\sigma)}=\pi_{i,j}(N_\sigma).$

For  $\sigma\in {\bf S}_n^2$, the matrix $N_{\sigma}\in M_{n\times
n}(\mathbb{Z}^+)$ defined (\ref{eq1}), can be visualized as a set of
integral points of the corresponding coordinates in the lattice of
integral points of $\mathbb{R}^2$. Respectively we can visualize
$\sigma$ by considering $\sigma=(i_1,j_1)\ldots(i_k,j_k)$ where
$i_s<j_s$ for any $s\ :\ 1\leq s\leq k$ as a set of $k$ integral
points with coordinates $(i_s,j_s).$ By slight abuse of notation we
will not distinguish between $\sigma$ and its visualization.

By the definition of a rank matrix, $(R_\sigma)_{(i,j)}$ is the
number of points of $\sigma$ inside the rectangular triangle with
vertices $(i,i)$, $(j,j)$ and $(i,j)$ (including the points on the
legs of the triangular).

Now, given subset $P\subseteq \sigma \subseteq \mathbb{R}^2$ let
$\langle P\rangle$ to be its support. Set $T(P)$ to be the smallest
rectangular triangle in $\mathbb{R}^2$ containing $P$ such that its
hypotenuse lies on the line $y=x$. Set also
$T(\emptyset):=\emptyset$. Obviously $T(P)$ is completely determined
by its rectangular vertex which has coordinates $(x_P,y_P)$ where
$x_P=\min\{i;\ (i,j)\in P\}=\min \langle P\rangle$ and
$y_P=\max\{j;\ (i,j)\in P\}=\max \langle P\rangle$. Let $\hat P$ be
the set of all points of $\sigma$ in $T(P).$ Then
$$
{\rm card}(\hat
P)=(R_{\sigma})_{(x_{P},y_{P})}=L(\pi_{x_P,y_P}(\sigma))
$$

and the order $\preceq$ can be immediately translated by

\begin{lem}\label{lem1}
The following claims are equivalent
\begin{enumerate}
\item[(i)]
$\sigma'\preceq \sigma$.
\item [(ii)] for any subset $P'\subseteq
\sigma'$ there exists a subset $P\subseteq \sigma$ such that
$T(P)\subseteq T(P')$ and ${\rm card} (\hat P)\geq {\rm card} (\hat
P')$.
\end{enumerate}
\end{lem}

\begin{proof}
(i)$\Rightarrow$(ii) For any subset $P'\subseteq \sigma'$, we have
${\rm card}(\hat P')=(R_{\sigma'})_{(x_{P'},y_{P'})}$. Then it is
enough to consider $P:=T(P')\cap \sigma$.

(ii)$\Rightarrow$(i) For any integers $1\leq i<j \leq n$, let
$T_{(i,j)}$ be a rectangular triangular with vertices $(i,j),\
(i,i),\ (j,j)$. Consider $P':=\sigma'\cap T$ (thus, $\hat P'=P'$).
Then by hypothesis there exists a subset $P\subseteq \sigma$ such
that $T(P)\subseteq T(P')$ and ${\rm card} (\hat P)\geq {\rm card} (
P')$. In particular we have $\hat P\subseteq T_{(i,j)}\cap \sigma$,
and therefore $(R_{\sigma})_{i,j}={\rm card}(T_{(i,j)}\cap
\sigma)\geq {\rm card} (\hat P)\geq{\rm card}
(P')=(R_{\sigma'})_{i,j}$.
\end{proof}

\subsection{}\label{4.6b}
Consider an involution $\sigma=(i_1,j_1)\ldots(i_k,j_k).$  There are
two kinds of {\bf elementary transformations} giving a new
involution in $\bS_n^2(k):$

\begin{enumerate}
\item[(I)] We change an integer\ \ $p\in \langle \sigma\rangle$  by an integer
$q\in \langle \sigma\rangle^c$. Geometrically its means that we move
a point $(i,j)$ of $\sigma$ to the point $(i',j)$ where $i'\in
\langle \sigma\rangle^c$ (that is vertically) or to the point
$(i,j')$ where $j'\in \langle \sigma\rangle^c$ (that is
horizontally). We denote this transformation by $F_{(p\rightarrow
q)}(\sigma).$

\item[(II)] We exchange two integers $p,q\in\langle \sigma\rangle$ from different
pairs of $\sigma$:
  \begin{enumerate}
    \item[a)-] If  $p=j_s$ and $q=i_t$, geometrically
    its means that we move one point vertically and the other point
    horizontally.

    \item[b)-] If  $p=i_s$ and $q=i_t$, geometrically
    its means that we move vertically the two corresponding points.
     [Notice the
exchange of $j_s$ and $j_t$ is exactly the same operation, therefore
we can see horizontally the movement of the two points].
    \end{enumerate}
    We denote this transformation by $F_{(p\rightleftarrows q)}(\sigma)$
\end{enumerate}

Of course the direction of the displacement of the point (or of the
two points) is important, the resulting involution for one
elementary transformation will be smaller or bigger than $\sigma$
for the order $\preceq$. Namely
\begin{enumerate}
\item [(I)] If the elementary transformation is of type (I)
 then moving  a point\ up or to the right increases the sizes of the
 rectangular triangles for
the resulting involution, so that by Lemma \ref{lem1} (ii) the
resulting involution is smaller. Respectively moving a point down or
to the left gives a bigger involution. In other words
$$
F_{(i_s\rar x)}(\sigma)\left\{\begin{array}{ll}
                           \prec\sigma & {\rm if\ }x<i_s,\\
                           \succ\sigma &{\rm otherwise.}\\
                           \end{array}
                           \right.\quad{\rm and}\quad
F_{(j_s\rar y)}(\sigma)\left\{\begin{array}{ll}
                           \prec\sigma & {\rm if\ }y>j_s,\\
                           \succ\sigma &{\rm otherwise.}\\
                           \end{array}\right.
\eqno{(*)}
$$

\item [(II)] If the elementary transformation is of type (II) we  put $i_s<i_t$
and use the figures to comment the situation. In figures below we
denote by $\bullet$ the points of $\sigma$ and by $\times$ the
points of $F(\sigma)$.
     \begin{enumerate}

\item [a)-] We exchange $j_s$ and $i_t.$ If $j_s<i_t$ (resp.
        $j_s<i_t$) then  the distribution of points
        will increase the sizes of the different
rectangular triangles and therefore by Lemma \ref{lem1} (ii) the
resulting involution will be smaller (resp. bigger),
 [see Figure \ref{Fig1}],

\begin{figure}[htbp]
  \centering
  \psset{unit=0.5cm}
  \begin{tabular}{cc}

  \begin{pspicture}(-15,-1)(13,3)

\psline(-9.6,-0.5)(-6.6,2.5)

 \put(-6,2.5){$(i_t,j_t)$}
\psline[linestyle=dotted,dotsep=1pt](-4.5,2)(-7,2)(-7,-2)
\psline[linestyle=dotted,dotsep=1pt](-4.5,0.5)(-8.5,0.5)(-8.5,-2)

 \put(-8.7,-1){$\bullet$}
 \put(-12.5,-1.5){$(i_s,j_s)$}
 \put(-7.3,-1){$\times$}

\put(-5.6,0.2){$\times$}
 \put(-5.4,1.8){$\bullet$}
 \pscurve[linestyle=dashed,
 linearc=1pt]{->}(-4.8,1.8)(-4.3,1.15)(-4.8,0.5)

\pscurve[linestyle=dashed,
  linearc=1pt]{->}(-8.7,-0.8)(-7.9,-1.6)(-7.1,-0.8)

  \rput{U}(-0.5,-1){resp.}

\psline(3.1,-0.5)(6.1,2.5)

\psline[linestyle=dotted,dotsep=1pt](8.1,2)(5.6,2)(5.6,-2)
\psline[linestyle=dotted,dotsep=1pt](8.1,0.5)(4.1,0.5)(4.1,-2)

 \put(3.8,-1){$\times$}
 \put(5.4,-1){$\bullet$}
 \put(7.5,-2){$(i_s,j_s)$}

\put(7.2,0.2){$\bullet$}
 \put(7,1.8){$\times$}
 \put(8.8,0){$(i_t,j_t)$}
 \pscurve[linestyle=dashed,
 linearc=1pt]{<-}(7.8,1.8)(8.1,1.15)(7.8,0.5)

\pscurve[linestyle=dashed,
  linearc=1pt]{<-}(3.9,-0.8)(4.7,-1.6)(5.5,-0.8)

  \end{pspicture}
  \end{tabular}\caption{} \label{Fig1}
\end{figure}

\item[b)] Suppose that we exchange $i_s$ and $i_t$
         with $j_t<j_s$ (resp. $j_t>j_s$), again
the resulting involution will be smaller (resp. bigger) [see Figure
\ref{Fig2}].

\begin{figure}[htbp]
  \centering
\psset{unit=0.5cm}
\begin{tabular}{cc}
  \begin{pspicture}(-1,-1)(17,4)

\psline[linestyle=dotted,dotsep=1pt](0,0)(5,0)(5,3)(0,3)(0,0)
\put(-0.3,2.8){$\bullet$}
 \rput{U}(0,4){$(i_t,j_t)$}
 \put(4.7,-0.2){$\bullet$}
 \rput{U}(5,-1){$(i_s,j_s)$}
  \rput{U}(0,0){$\times$}
  \rput{U}(5,3){$\times$}
  \pscurve[linestyle=dashed,
  linearc=5pt]{<-}(-0.25,0.1)(-0.75,1.4)(-0.25,2.7)
  \pscurve[linestyle=dashed,
  linearc=1pt]{<-}(5.25,2.9)(5.75,1.4)(5.25,0.3)

  \rput{U}(8.5,1){resp.}

  \psline[linestyle=dotted,dotsep=1pt](12,0)(17,0)(17,3)(12,3)(12,0)
  \put(11.7,-0.2){$\bullet$} \put(16.65,2.8){$\bullet$}

  \rput{U}(17,0){$\times$}
  \rput{U}(17,4){$(i_t,j_t)$}

  \rput{U}(12,3){$\times$}
  \rput{U}(12,-1){$(i_s,j_s)$}

  \pscurve[linestyle=dashed,
  linearc=1pt]{<-}(11.75,2.9)(11.25,1.4)(11.75,0.3)

  \pscurve[linestyle=dashed,
  linearc=1pt]{<-}(17.25,0.1)(17.75,1.4)(17.25,2.7)

\end{pspicture}
\end{tabular}\caption{} \label{Fig2}
\end{figure}

     \end{enumerate}
\end{enumerate}

Note that by these elementary transformations we can pass from any
$\sigma\in\bS_n^2(k)$ to any other $\sigma'\in\bS_n^2(k).$

Moreover by \cite{Mel1} one gets

\begin{cor}\label{cor1} Let $\sigma,\sigma'\in \bS_n^2(k)$ be such
that $\sigma'\prec \sigma$  then there exists a sequence of
involutions $\sigma_0=\sigma,$ $\sigma_m=\sigma'$ such that for any
$i\ :\ 1\leq i\leq m$ $\sigma_i$ is obtained from $\sigma_{i-1}$ by
an elementary transformation and  $\sigma_0\succ
\sigma_1\succ\sigma_2\succ\cdots\succ\sigma_m$.
\end{cor}

By this corollary,  any descendant (or ancestor) of $\sigma$ is
obtained from $\sigma$ by an elementary transformation which cannot
be decomposed into a sequence of other elementary transformations
preserving the order $\preceq$. We will call it a {\bf minimal}
elementary transformation.

\subsection{}\label{4.61} In the next four  subsections we describe 4
types of minimal elementary transformations. All the proofs are
given in {Mel1}, here we only give some explanations.

Let us first consider the elementary transformation of type (I)
moving a point vertically and giving smaller involution. Thus, we
move a point $P:(i_s,j_s)$ of $\sigma$ down to
 $P':(x,j_s)$ where $x\in \langle \sigma\rangle^c$ [i.e. the
elementary transformation is $F_{(i_s\rar x)}(\sigma)$]. If there
exists $y\in \langle \sigma\rangle^c$ such that $i_s>y>x$ then
$F_{(i_s\rar x)}$ can be decomposed into 2 elementary
transformations as one can see from Fig. \ref{Fig3a}.

\begin{figure}[htbp]
  \centering
\psset{unit=0.5cm}
\begin{tabular}{cc}
  \begin{pspicture}(-1,-1)(20,4)

 \put(0,0){$\times$}  \put(1,0){$ (x,j_s)$}
 \put(0,1.5){$\times$}\put(1, 1.5){$(y,j_s)$}
 \put(0,3){$\bullet$} \put(1,3){$(i_s, j_s)$}
\pscurve[linestyle=dashed,linearc=1pt]{<-}(-0.15,0.35)(-0.4,1.65)(-0.15,2.8)

\put(5.5,1.5){=}

 \put(8,0){$\times$}
 \put(8,1.5){$\times$}
 \put(8,3){$\bullet$}
\pscurve[linestyle=dashed,linearc=1pt]{->}(7.75,2.8)(7.5,2.2)(7.75,1.6)

\put(10,1.5){$\oplus$}

\put(12,0){$\times$}
 \put(12,1.5){$\bullet$}
 \put(12,3){$\times$}

\pscurve[linestyle=dashed,linearc=1pt]{<-}(11.75,0.6)(11.5,1.2)(11.75,1.8)
\end{pspicture}
\end{tabular}
\caption{} \label{Fig3a}
\end{figure}

Thus, the necessary condition for $F_{(i_s\rar x)}$ to be minimal is
$x=\max\{p\in \langle \sigma\rangle^c\ :\ p<i_s\}.$ However even in
this case $F_{(i_s\rar x)}$ may be not minimal. Indeed if there
exists $(i_t,j_t)\in\sigma$ such that $x<i_t<i_s$ and $j_t>j_s$ then
 $F_{(i_s\rar x)}$ can be
decomposed into 3 elementary transformations as one can see from
Fig. \ref{Fig3}.

\begin{figure}[htbp]
  \centering
\psset{unit=0.5cm}
\begin{tabular}{cc}
  \begin{pspicture}(-1,-1)(20,4)

\put(0,0){$\times$}  \put(0,-0.75){$(x,j_s)$}
\put(-0.1,2.75){$\bullet$} \put(0,3.75){$(i_s,j_s)$}
\pscurve[linestyle=dashed,linearc=1pt]{<-}(-0.15,0.35)(-0.4,1.65)(-0.15,2.8)

\put(1.5,1.5){$\bullet$}  \put(2,1.5){$(i_t,j_t)$}

\put(4.5,1.2){=}

 \put(6.75,1.5){$\times$}
 \put(6.8,2.75){$\bullet$}
\pscurve[linestyle=dashed,linearc=1pt]{<-}(6.6,1.6)(6.2,2.1)(6.6,2.65)
 \put(8.45, 2.75){$\times$}
\put(8.5,1.5){$\bullet$}\put(6.75,0){$\times$}
\pscurve[linestyle=dashed,linearc=1pt]{->}(9,1.6)(9.5,2.1)(9,2.65)

\put(10.5,1.2){$\oplus$}

\put(12.5,1.5){$\bullet$}  \put(12.4,0){$\times$}
\pscurve[linestyle=dashed,linearc=1pt]{->}(12.25,1.65)(11.8,0.85)(12.25,0.2)

\put(14,2.75){$\bullet$}

\put(15.5,1.5){$\oplus$}

\put(17.5,0){$\bullet$}

\put(19.5,3){$\bullet$} \put(19.4,1.5){$\times$}
\pscurve[linestyle=dashed,linearc=1pt]{->}(20.1,2.8)(20.5,2.2)(20.1,1.6)

\end{pspicture}
\end{tabular}\caption{} \label{Fig3}
\end{figure}
So the second necessary condition is that for any
$(i_t,j_t)\in\sigma$ such that $x<i_t<i_s$ one has $j_t<j_s.$
Moreover, by \cite[3.11]{Mel1} these two conditions are sufficient
for the existence of  the minimal elementary vertical transformation
giving a smaller involution.

Formally, the minimal elementary vertical transformation giving a
smaller involution exists if
\begin{itemize}
\item{} $x:=\max\{\langle
\sigma\rangle^c\cap\{i\}_{i=1}^{i_s-1}\} $ exists, (i.e.
  $\langle
\sigma\rangle^c\cap\{i\}_{i=1}^{i_s-1}\ne\emptyset$);
 \item{} Either $x=i_s-1$ or for any point $(i_t,j_t)\in \sigma$ such that
$x<i_t<i_s$ one has $j_t<j_s.$
\end{itemize}
In that case put $\sigma_{i_s\dar}:=F_{i_s\rar x}(\sigma).$
Otherwise put $\sigma_{i_s\dar}:=\emptyset.$
\begin{rem}\label{remark}
Note that by our definition if  $\sigma_{i_s\dar}=F_{i_s\rar
x}(\sigma)$ then  for any $t: x<t<i_s$ one has either $t=j_p$ for
some pair $(i_p,j_p)\in \sigma$ or $t=i_p$ for some pair
$(i_p,j_p)\in\sigma$ such that $j_p<j_s.$
\end{rem}

Respectively, the minimal elementary transformation moving
$(i_s,j_s)$ up giving a bigger involution exists if
\begin{itemize}
\item{} $x:=\min\{\langle
\sigma\rangle^c\cap\{i\}_{i=i_s+1}^{j_s-1}\} $ exists, (i.e.
  $\langle
\sigma\rangle^c\cap\{i\}_{i=i_s+1}^{j_s-1}\ne\emptyset$);
 \item{} Either $x=i_s+1$ or for any point $(i_t,j_t)\in \sigma$ such that
$i_s<i_t<x$ one has $j_t<j_s.$
\end{itemize}
In that case put $\sigma_{\uar i_s}:=F_{i_s\rar x}(\sigma).$
Otherwise put $\sigma_{\uar i_s}:=\emptyset.$

Note that the two minimal elementary transformations described above
are inverse to each other:  if $\sigma_{\uar i_s}=F_{(i_s\rar
x)}(\sigma)$ then $(\sigma_{\uar i_s})_{x\dar}=\sigma$ and if
$\sigma_{i_s\dar}=F_{(i_s\rar x)}(\sigma)$ then
$(\sigma_{i_s\dar})_{\uar x}=\sigma.$

\begin{ex} Take $\sigma=(2,6)(3,5)(7,9)(8,10)\in\bS_{11}^2.$ Then

\noindent
 $\sigma_{2\dar}=F_{2\rar
1}(\sigma)=(1,6)(3,5)(7,9)(8,10),$

\noindent
 $\sigma_{\uar 2 }=F_{2\rar
4}(\sigma)=(4,6)(3,5)(7,9)(8,10),$\ \ \ $\sigma_{3\dar}=\emptyset,$\

\noindent
 $\sigma_{\uar 3}=F_{3\rar
4}(\sigma)=(2,6)(4,5)(7,9)(8,10),$\

\noindent
 $\sigma_{7\dar}=F_{7\rar
4}(\sigma)=(2,6)(3,5)(4,9)(8,10),$\ \ $\sigma_{\uar 7}=\emptyset,$ \

\noindent
 $\sigma_{8\dar}=F_{8\rar
4}(\sigma)=(2,5)(3,4)(7,9)(4,10),$ \ \  $\sigma_{\uar 8
}=\emptyset.$
\end{ex}

\subsection{}\label{4.62}
Considering a symmetry around $y=-x+n$ we get a minimal elementary
transformation moving $(i_s,j_s)$ horizontally, that is $F_{(j_s\rar
y)}(\sigma).$

By \cite[3.12]{Mel1}  (or by the symmetry mentioned above) the
minimal elementary transformation moving a $(i_s,j_s)$ to the right
giving a smaller involution exists if
 \begin{itemize}
 \item{} $y:=\min\{\langle \sigma\rangle^c\cap\{j\}_{j=j_s+1}^{n}\}$
 exists (i.e $\langle
\sigma\rangle^c\cap\{j\}_{j=j_s+1}^{n}\ne\emptyset$);
 \item{}  Either $y=j_s+1$ or for any point $(i_t,j_t)\in \sigma$ such that
$j_s<j_t<y$ one has $i_t>i_s.$
\end{itemize}
 In that case put $\sigma_{j_s\rar}:=F_{(j_s\rar y)}(\sigma).$
 Otherwise put $\sigma_{j_s\rar}:=\emptyset.$

 Respectively, the minimal elementary transformation moving
$(i_s,j_s)$ to the left giving a bigger involution exists if
\begin{itemize}
 \item{} $y:=\max\{\langle \sigma\rangle^c\cap\{j\}_{j=i_s+1}^{j_s-1}\}$
 exists (i.e $\langle
\sigma\rangle^c\cap\{j\}_{j=i_s+1}^{j_s-1}\ne\emptyset$);
 \item{}  Either $y=j_s-1$ or for any point $(i_t,j_t)\in \sigma$ such that
$y<j_t<j_s$ one has $i_t>i_s.$
\end{itemize}
 In that case put $\sigma_{\lar j_s}:=F_{(j_s\rar y)}(\sigma).$
 Otherwise put $\sigma_{\lar j_s}:=\emptyset.$

\begin{rem}
Again that the two minimal elementary transformations described
above are inverse to each other: if $\sigma_{j_s\rar}=F_{(j_s\rar
y)}(\sigma)$ then $(\sigma_{j_s\rar})_{\lar y}=\sigma$ and if
$\sigma_{\lar j_s}=F_{(j_s\rar y)}(\sigma)$ then $(\sigma_{\lar
j_s})_{ y\rar}=\sigma.$
\end{rem}
\begin{ex} Take again $\sigma=(2,6)(3,5)(7,9)(8,10)\in \bS_{11}^2.$
Then

\noindent $\sigma_{6\rar}=F_{6\rar
11}(\sigma)=(2,11)(3,5)(7,9)(8,10),$\

\noindent
 $\sigma_{\lar 6}=F_{6\rar
4}(\sigma)=(2,4)(3,5)(7,9)(8,10),$\ \ $\sigma_{5\rar}=\emptyset,$\

\noindent
 $\sigma_{\lar 5}=F_{5\rar
4}(\sigma)=(2,6)(3,4)(7,9)(8,10),$\

\noindent
 $\sigma_{9\rar}=F_{9\rar 11}(\sigma)
=(2,6)(3,5)(7,11)(8,10),$\ \ $\sigma_{\lar 9}=\emptyset,$\

\noindent
 $\sigma_{10\rar}=F_{10\rar
11}(\sigma)=(2,6)(3,5)(7,9)(8,11)$ \ \  $\sigma_{\lar
10}=\emptyset.$
\end{ex}

\subsection{}\label{4.63}
Now we turn to the minimal elementary transformations of type
(II)-a) giving a smaller involution. Let $(i_s,j_s)$ and $(i_t,j_t)$
be two points of $\sigma$ such that $i_s<i_t$. Then by Fig.
\ref{Fig1} $F_{(j_s\rightleftarrows i_t)}(\sigma)<\sigma$ if
$j_s<i_t.$

If $\langle \sigma\rangle^c\cap\{r\}_{j_s+1}^{i_t-1}\ne\emptyset$
then $F_{(j_s\rightleftarrows i_t)}(\sigma)<\sigma$ is not minimal.
Indeed let $y\in \langle \sigma\rangle^c\cap\{r\}_{j_s+1}^{i_t-1}$
then  $F_{(j_s\rightleftarrows i_t)}(\sigma)$ can be decomposed into
3 elementary transformations as one can see from Fig. \ref{Fig4a}.
\begin{figure}[htbp]
 \centering
  \psset{unit=0.5cm}
  \begin{pspicture}(-5,-4)(13,3)

\psline(-9.6,-0.5)(-6.6,2.5)

 \put(-6,2.5){$(i_t,j_t)$}
\psline[linestyle=dotted,dotsep=1pt](-4.5,2)(-7,2)(-7,-2)
\psline[linestyle=dotted,dotsep=1pt](-4.5,0.5)(-8.5,0.5)(-8.5,-2)
\psline[linestyle=dotted,dotsep=1pt](-4.5,1.25)(-7.75,1.25)(-7.75,-2)
 \put(-8.7,-2.5){$j_s$}
 \put(-7,-2.5){$i_t$}
 \put(-8,-2.5){$y$}

 \put(-8.7,-1){$\bullet$}
 \put(-11.2,-1.5){$(i_s,j_s)$}
 \put(-8.1,-1){$\times$}
 \put(-7.3,-1){$\times$}

\put(-5.6,0.2){$\times$}
 \put(-5.4,1.8){$\bullet$}
 \pscurve[linestyle=dashed,
 linearc=1pt]{->}(-4.8,1.8)(-4.3,1.15)(-4.8,0.5)

\pscurve[linestyle=dashed,
  linearc=1pt]{->}(-8.7,-0.8)(-7.9,-1.6)(-7.1,-0.8)
 \put(-3,1){$=$}

\psline(-1.5,-0.5)(1.5,2.5)

\put(2.5,1.8){$\bullet$}
\psline[linestyle=dotted,dotsep=1pt](3.5,2)(1,2)(1,-2)
\psline[linestyle=dotted,dotsep=1pt](3.5,0.5)(-0.5,0.5)(-0.5,-2)
\psline[linestyle=dotted,dotsep=1pt](3.5,1.25)(0.25,1.25)(0.25,-2)

 \put(-0.4,-2.5){$j_s$}
 \put(1,-2.5){$i_t$}
 \put(0.2,-2.5){$y$}

 \put(-0.7,-1){$\bullet$}
 \put(0,-1){$\times$}
 \put(0.8,-1){$\times$}

\pscurve[linestyle=dashed,
 linearc=1pt]{->}(-0.5,-1.2)(-0.05,-1.6)(0.4,-1.2)

  \put(4.5,0.9){$\oplus$}
\psline(5,-0.5)(8,2.5)

\psline[linestyle=dotted,dotsep=1pt](10,2)(7.5,2)(7.5,-2)
\psline[linestyle=dotted,dotsep=1pt](10,1.25)(6.75,1.25)(6.75,-2)

 \put(7.5,-2.5){$i_t$}
 \put(6.6,-2.5){$y$}

 \put(6.6,-1){$\bullet$}
 \put(7.2,-1){$\times$}

\put(9,1.8){$\bullet$}
 \put(8.9,1){$\times$}

 \pscurve[linestyle=dashed,
 linearc=1pt]{->}(6.7,-1.2)(7.1,-1.6)(7.5,-1.2)

\pscurve[linestyle=dashed,
 linearc=1pt]{->}(9.5,2)(9.8,1.6)(9.5,1.2)
  \put(11,0.9){$\oplus$}
\psline(11.5,-0.5)(14.5,2.5)
\psline[linestyle=dotted,dotsep=1pt](16,2)(14,2)(14,-2)
\psline[linestyle=dotted,dotsep=1pt](16.5,0.5)(12.5,0.5)(12.5,-2)
\psline[linestyle=dotted,dotsep=1pt](16.5,1.25)(13.25,1.25)(13.25,-2)

\put(15,1.1){$\bullet$}
 \put(14.9,0.3){$\times$}
 \put(13.8,-1){$\bullet$}
 \put(14.2,-1.5){$(i_s,i_t)$}

\put(13.2,-2.5){$y$}
 \put(12.3,-2.5){$j_s$}
  \put(14,-2.5){$i_t$}
\put (16.7,0.3){$(j_s,j_t)$}
 \pscurve[linestyle=dashed,
 linearc=1pt]{->}(15.3,1.2)(15.7,0.8)(15.3,0.4)
\end{pspicture}
\caption{} \label{Fig4a}
\end{figure}

Thus, the necessary condition for $F_{(j_s\rightleftarrows
i_t)}(\sigma)$ to be minimal is $\langle
\sigma\rangle^c\cap\{r\}_{j_s+1}^{i_t-1}=\emptyset.$ Moreover,
assume that there exists a point $(i,j)$ of $\sigma$ such that
either $j_s<i<i_t$ and $j>j_s$ or $j_s<j<i_t$ and $i<i_s$ (i.e. a
point $(i,j)$ outside of the triangular with vertices $(i_s,i_s),\
(j_t,j_t),\ (i_s,j_t)$ such that either $x=i$ or $y=j$ crosses the
sides of the triangular). Then $F_{(j_s\rightleftarrows
i_t)}(\sigma)$ can be decomposed into 3 elementary transformations.
Since these cases are obtained one from another by symmetry around
$y=n-x$ we show in Fig. \ref{Fig4} only the case $j_s<i<i_t$ and
$j>j_t.$

\begin{figure}[htbp]
  \centering
\psset{unit=0.6cm}
  \begin{pspicture}(-12,-4)(13,5)
  \psline(-10.7,-0.5)(-7.2,3)

 \put(-5,4){$(i_t,j_t)$}
\psline[linestyle=dotted,dotsep=1pt](-5,2.5)(-7.7,2.5)(-7.7,-2.2)
\psline[linestyle=dotted,dotsep=1pt](-5.1,0.8)(-9.4,0.8)(-9.4,-2.2)

 \put(-9.3,-2.5){$j_s$}
 \put(-7.6,-2.5){$i_t$}

 \put(-9.5,-0.8){$\bullet$}
 \put(-11.8,-2){$(i_s,j_s)$}
 \put(-8,-0.8){$\times$}

\put(-6.2,0.6){$\times$}
 \put(-6,2.4){$\bullet$}
\put(-3,1.5){$\bullet$}
 \put(-3,2.2){$(i,j)$}
 \pscurve[linestyle=dashed,
 linearc=1pt]{->}(-5.8,2.3)(-5.3,1.6)(-5.8,0.9)

\pscurve[linestyle=dashed,
  linearc=1pt]{->}(-9.2,-0.9)(-8.6,-1.3)(-8,-0.9)
 \put(-1,1){$=$}
 \psline(-1.1,-0.5)(2.6,3)

 \put(2,4){$(i_t,j_t)$}
%
\psline[linestyle=dotted,dotsep=1pt](4.6,0.8)(0.3,0.8)(0.3,-2.2)


 \put(0.2,-0.8){$\bullet$}
 \put(-1.5,-2){$(i_s,j_s)$}

 \put(3.6,2.4){$\bullet$}
 \put(3.5,1.5){$\times$}
\put(5.5,1.5){$\bullet$}
 \put(5.4,2.4){$\times$}
 \put(4.5,-0.5){$(i,j)$}
 \pscurve[linestyle=dashed,
linearc=1pt]{->}(3.4,2.3)(3,1.95)(3.4,1.6)

\pscurve[linestyle=dashed,
 linearc=1pt]{->}(6,1.6)(6.4,1.95)(6,2.3)
\put(8,1){$\oplus$}

\end{pspicture}
\end{figure}
\begin{figure}[htbp]
 \centering
  \psset{unit=0.6cm}
  \begin{pspicture}(-12,-4)(13,3)

\psline(-10.8,-0.5)(-7.3,3)
 \put(-4.7,3.5){$(i_t,j)$}
\psline[linestyle=dotted,dotsep=1pt](-4.7,1.8)(-8.4,1.8)(-8.4,-1.5)
\psline[linestyle=dotted,dotsep=1pt](-4.7,0.8)(-9.4,0.8)(-9.4,-1.5)


 \put(-9.5,-0.4){$\bullet$}
 \put(-8.7,-0.4){$\times$}
 \put(-9.5,-2.5){$j_s$}
  \put(-8.5,-2.5){$i$}
 \put(-12.2,-1){$(i_s,j_s)$}

 \put(-6.7,1.6){$\bullet$}
 \put(-6.9,0.5){$\times$}
 \put(-8.2,3){$(i,j_t)$}
 \put(-4.2,2.4){$\bullet$}
 \pscurve[linestyle=dashed,
linearc=1pt]{->}(-6.3,1.6)(-5.9,1.15)(-6.3,0.7)

\pscurve[linestyle=dashed,
 linearc=1pt]{->}(-9.4,-1)(-9,-1.4)(-8.6,-1)
\put(-1,2){$\oplus$}
 \psline(-1.1,-0.5)(2.6,3)
\psline[linestyle=dotted,dotsep=1pt](6.2,2.5)(2,2.5)(2,-1.6)
\psline[linestyle=dotted,dotsep=1pt](6.2,1.5)(1,1.5)(1,-1.6)

 \put(0.9,-2.5){$i$}
 \put(2.1,-2.5){$i_t$}

 \put(0.9,-0.3){$\bullet$}
 \put(1.8,-0.3){$\times$}
 \put(-1.3,-0.5){$(i_s,i)$}

\put(3.5,0.6){$\bullet$}
 \put(3.7,-0.5){$(j_s,j_t)$}
\put(5.2,2.4){$\bullet$}
 \put(5.2,1.3){$\times$}
 \put(6.7,3){$(i_t,j)$}
 \pscurve[linestyle=dashed,
 linearc=1pt]{->}(5.5,2.4)(5.9,1.9)(5.5,1.3)

\pscurve[linestyle=dashed,
  linearc=1pt]{->}(1,-0.5)(1.35,-0.9)(1.8,-0.5)

%
\end{pspicture}
\caption{} \label{Fig4}
\end{figure}
Thus the necessary conditions for $F_{(j_s\rightleftarrows
i_t)}(\sigma)$ to be minimal are
\begin{itemize}
\item{} $\langle \sigma\rangle^c\cap\{p\}_{j_s+1}^{i_t-1}=\emptyset$
\item{} either $j_s=i_t-1$ or for
any $q:\ j_s<q<i_t$ one has $q\in
\langle\pi_{i_s,j_t}(\sigma)\rangle.$
\end{itemize}
By \cite[3.13]{Mel1} these conditions are also sufficient, that is
if $(i_s,j_s),\ (i_t,j_t)$ satisfy the conditions above then
$F_{(j_s\rightleftarrows i_t)}(\sigma)$ is minimal.

Summarizing, to get all the minimal elementary transformations of
type (II)-a) giving a smaller element for a given
$(i_t,j_t)\in\sigma$ put $E_{(i_t,j_t)}(\sigma)$ be the set of all
pairs $(i_s,j_s)\in\sigma$ such that $(i_s,j_s),\ (i_t,j_t)$ satisfy
the two conditions above (note, that there can be a few such pairs
$(i_s,j_s)$).

 If $E_{(i_t,j_t)}\ne\emptyset$ put
$C_{i_t\looparrowright}(\sigma):=\{F_{j_s\rightleftarrows
i_t}(\sigma)\ :\
 (i_s,j_s)\in E_{(i_t,j_t)}(\sigma)\}.$ If
 $E_{(i_t,j_t)}(\sigma)=\emptyset$ put $C_{i_t\looparrowright}(\sigma):=\emptyset.$

\begin{rem}\label{rem1}
Note that by the conditions above if $F_{j_s\rightleftarrows
i_t}(\sigma)$ is minimal then
\begin{itemize}
\item[(i)] For any $(i_p,j_p)\in\sigma$ such that $i_p<i_s$ one has either $j_p<j_s$ or
$j_p>i_t.$
\item[(ii)] For any $(i_p,j_p)\in\sigma$ such that $i_s<i_p<i_t$ one has
either $i_p<j_s$ or $j_p<j_t.$
\item[(iii)] If $i_t\ne j_s+1$ then for any $x\in \{a\}_{a=j_s+1}^{i_t-1}$
one has $x\in\langle\sigma\rangle$ and for $(i_p,j_p)\in\sigma$ such
that $x\in\{i_p,j_p\}$ one has $i_s<i_p<j_p<j_t.$
\end{itemize}
\end{rem}

Respectively, to get all the minimal elementary transformations of
type (II)-a) giving a bigger element for a given
$(i_t,j_t)\in\sigma$ put $G_{(i_t,j_t)}(\sigma)$ to be the set of
all $(i_s,j_s)\in\sigma$ such that $i_s<i_t$,\ $j_s>i_t$ and
\begin{itemize}
\item{} $\langle \sigma\rangle^c\cap\{p\}_{i_t+1}^{j_s-1}=\emptyset$
\item{} either $j_s=i_t+1$ or for
any $q:\ i_t<q<j_s$ one has $q\in
\langle\pi_{i_t,j_s}(\sigma)\rangle.$
\end{itemize}

 If
$G_{(i_t,j_t)}\ne\emptyset$ put $C_{\looparrowleft
i_t}(\sigma):=\{F_{i_t\rightleftarrows j_s}(\sigma)\ :\ (i_s,j_s)\in
G_{(i_t,j_t)}(\sigma)\}.$ If
 $G_{(i_t,j_t)}(\sigma)=\emptyset$ put $C_{\looparrowleft i_t}(\sigma):=\emptyset.$

\begin{ex} Take $\sigma=(1,3)(2,4)(5,9)(6,10)(7,8).$ Then
$$\begin{array}{lcl}
C_{2\looparrowright}(\sigma)&=&\emptyset,\\
 C_{\looparrowleft
2}&=\{&F_{(2\rightleftarrows 3)}(\sigma)=(1,2)(3,4)(5,9)(6,10)(7,8)\ \},\\
C_{5\looparrowright}(\sigma)&=\{&F_{(4\rightleftarrows
5)}(\sigma)=(1,3)(2,5)(4,9)(6,10)(7,8),\\
&&F_{(3\rightleftarrows 5)}(\sigma)=(1,5)(2,4)(3,9)(6,10)(7,8)\ \},\\
C_{\looparrowleft 5}(\sigma)&=&\emptyset\\
C_{6\looparrowright}(\sigma)&=\{&F_{(4\rightleftarrows
6)}(\sigma)=(1,3)(2,6)(5,9)(4,10)(7,8),\\
&&F_{(3\rightleftarrows 6)}(\sigma)=(1,6)(2,4)(5,9)(3,10)(7,8)\ \},\\
C_{\looparrowleft 6}(\sigma)&=\{&F_{(6\rightleftarrows 9)}(\sigma)=(1,3)(2,4)(5,6)(9,10)(7,8)\ \},\\
C_{7\looparrowright}(\sigma)&=&\emptyset,\\ C_{\looparrowleft 7}&=&\emptyset.\\
\end{array}$$
\end{ex}

Exactly as in previous cases the actions $\looparrowright$ and
$\looparrowleft$ are inverse, that is if
$E_{(i_t,j_t)}(\sigma)\ne\emptyset$ and $F_{(j_s\rightleftarrows
i_t)}(\sigma)\in C_{i_t\looparrowright}(\sigma)$ then
$F_{(j_s\rightleftarrows i_t)}(F_{(j_s\rightleftarrows
i_t)}(\sigma))=\sigma\in C_{\looparrowleft
j_s}(F_{(j_s\rightleftarrows i_t)}(\sigma));$ if
$G_{(i_t,j_t)}(\sigma)\ne \emptyset$ and $F_{(i_t\rightleftarrows
j_s)}(\sigma)\in C_{\looparrowleft i_t}(\sigma)$ then
$F_{(i_t\rightleftarrows j_s)}(F_{(i_t\rightleftarrows
j_s)}(\sigma))=\sigma\in
C_{j_s\looparrowright}(F_{(i_t\rightleftarrows j_s)}(\sigma).$

\subsection{}\label{4.64}
Finally consider  the minimal elementary transformations of type
(II)- b) giving a smaller involution. Let $(i_s,j_s)$ and
$(i_t,j_t)$ be two points of $\sigma$ such that $i_s<i_t$. Then by
Fig. \ref{Fig2} $F_{(i_s\rightleftarrows i_t)}(\sigma)<\sigma$ if
$j_s>j_t.$

Assume there exists $(i,j)\in\sigma$ such that $i_s<i<i_t$ and
$j_t<j<j_s.$ Then $F_{(i_s\rightleftarrows i_t)}(\sigma)$ can be
decomposed into 3 elementary transformations as one can see from
Fig. \ref{Fig5}.

\begin{figure}[htbp]
  \centering
\psset{unit=0.5cm}
  \begin{pspicture}(10,-1)(17,4)

\psline[linestyle=dotted,dotsep=1pt](-0.3,0)(3.7,0)(3.7,3)(-0.3,3)(-0.3,0)
\put(-0.6,2.8){$\bullet$}
 \put(1.2,1.3){$\bullet$}
 \put(1.7,1.3){$(i,j)$}
 \rput{U}(-0.3,4){$(i_t,j_t)$}
 \put(3.4,-0.2){$\bullet$}
 \rput{U}(3.7,-1){$(i_s,j_s)$}
  \rput{U}(-0.3,0){$\times$}
  \rput{U}(3.6,3){$\times$}
  \pscurve[linestyle=dashed,
  linearc=5pt]{<-}(-0.55,0.1)(-1.05,1.4)(-0.55,2.7)
  \pscurve[linestyle=dashed,
  linearc=1pt]{<-}(3.95,2.9)(4.45,1.4)(3.95,0.3)
\put(5.5,1.5){=}
\psline[linestyle=dotted,dotsep=1pt](8,3)(9.5,3)(9.5,1.3)(8,1.3)(8,3)
\put(7.7,2.8){$\bullet$}
 \put(7.7,1.2){$\times$}
 \put(9.3,1.2){$\bullet$}
 \put(9.2,2.8){$\times$}
 \put(9.5,0.3){$(i,j)$}
 \rput{U}(8,4){$(i_t,j_t)$}
 \put(11.5,-0.2){$\bullet$}
 \rput{U}(11.8,-1){$(i_s,j_s)$}
  \pscurve[linestyle=dashed,
  linearc=5pt]{<-}(7.4,1.6)(7.1,2.1)(7.4,2.8)
  \pscurve[linestyle=dashed,
  linearc=1pt]{->}(10,1.6)(10.3,2.1)(10,2.8)
\put(13,1.5){$\oplus$}
\psline[linestyle=dotted,dotsep=1pt](15,1.5)(19,1.5)(19,0)(15,0)(15,1.5)
 \put(14.8,1.3){$\bullet$}
 \put(14.7,-0.2){$\times$}
 \put(16.5,2.8){$\bullet$}
 \put(14,2){$(i,j_t)$}
 \rput{U}(16.5,4){$(i_t,j)$}
 \put(18.8,-0.2){$\bullet$}
 \put(18.7,1.3){$\times$}
 \rput{U}(19.3,-1){$(i_s,j_s)$}
  \pscurve[linestyle=dashed,
  linearc=5pt]{->}(14.6,1.2)(14.2,0.6)(14.6,0)
  \pscurve[linestyle=dashed,
  linearc=1pt]{->}(19.2,0)(19.6,0.6)(19.2,1.2)
\put(21,1.5){$\oplus$}
\psline[linestyle=dotted,dotsep=1pt](24.6,3)(27,3)(27,1.5)(24.6,1.5)(24.6,3)
 \put(23,-0.2){$\bullet$}
 \put(26.2,0.3){$(i,j_s)$}
 \rput{U}(25,4){$(i_t,j)$}
 \put(26.8,1.3){$\bullet$}
 \put(26.7,2.75){$\times$}
 \put(24.5,2.8){$\bullet$}
 \put(24.3,1.3){$\times$}
 \rput{U}(23,-1){$(i_s,j_t)$}
  \pscurve[linestyle=dashed,
  linearc=5pt]{->}(24.3,2.7)(23.9,2.1)(24.3,1.5)
  \pscurve[linestyle=dashed,
  linearc=1pt]{->}(27.2,1.5)(27.6,2.1)(27.2,2.7)

\end{pspicture}\caption{} \label{Fig5}
\end{figure}

Thus the necessary condition for $F_{i_s\rightleftarrows
i_t}(\sigma)$ to be minimal is that
\begin{itemize}
\item[$(*)$] For any
$(i,j)\in\sigma$ such that $i_s<i<i_t$ one has either $j<j_t$ or
$j>j_s$.
\end{itemize}
By \cite[3.14]{Mel1} this condition is also sufficient, that is
$F_{i_s\rightleftarrows i_t}(\sigma)$ is minimal iff the condition
is satisfied.

Summarizing, to get all the minimal elementary transformations of
type (II)-b) giving a smaller element for a given
$(i_s,j_s)\in\sigma$ set
$$K_{(i_s,j_s)}(\sigma):=\{(i_t,j_t) \ :\ i_t>i_s,\ j_t<j_s\ {\rm
and\ } i_s<i_q<i_t\ \Rightarrow\ j_q<j_t\ {\rm or}\ j_q>j_s \}$$
 If
$K_{(i_s,j_s)}(\sigma)\ne\emptyset$ put $C_{i_s\uar\dar
}(\sigma):=\{F_{(i_s\rightleftarrows i_t)}(\sigma)\ :\ (i_t,j_t)\in
K_{(i_s,j_s)}(\sigma)\}.$ If $K_{(i_s,j_s)}(\sigma)=\emptyset$ put
$C_{i_s\uar\dar}(\sigma):=\emptyset.$

Respectively, to get all the minimal elementary transformations of
type (II)-b) giving a bigger element for a given
$(i_s,j_s)\in\sigma$ set
$$
L_{(i_s,j_s)}(\sigma):=\{(i_t,j_t)\ :\ i_t>i_s,\ i_t<j_s<j_t\ {\rm
and}\  i_s<i_q<i_t\ \Rightarrow\ j_q<j_s\ {\rm or}\ j_q>j_t\}$$
 If $L_{(i_s,j_s)}(\sigma)\ne\emptyset$ put $C_{\dar\uar
 i_s}(\sigma):=\{F_{(i_s\rightleftarrows i_t)}(\sigma)\ :\ (i_t,j_t)\in
 L_{(i_s,j_s)}(\sigma).$ If $L_{(i_s,j_s)}(\sigma)=\emptyset$ put
$C_{\dar\uar i_s}(\sigma):=\emptyset.$

\begin{ex} Take $\sigma=(1,7)(2,5)(3,8)(4,6).$ Then
$$\begin{array}{l}
C_{1\uar\dar}(\sigma)= \{F_{(1\rightleftarrows
2)}(\sigma)=(1,5)(2,7)(3,8)(4,6),\
F_{(1\rightleftarrows 4)}(\sigma)=(1,6)(2,5)(3,8)(4,7)\},\\
C_{2\uar\dar}(\sigma)=\emptyset, \\
C_{3 \uar\dar}(\sigma)=\{F_{(3\rightleftarrows
4)}(\sigma)=(1,7)(2,5)(3,6)(4,8)\},\\
C_{4\uar\dar}(\sigma)=\emptyset;\\
C_{\dar\uar 1}=\{F_{(1\rightleftarrows 3)}(\sigma)=(1,8)(2,5)(3,7)(4,6)\}, \\
C_{\dar\uar 2}=\{F_{(2\rightleftarrows
3)}(\sigma)=(1,7)(2,8)(3,5)(4,6),\ F_{(2\rightleftarrows
4)}(\sigma)=(1,7)(2,6)(3,8)(4,5)\},\\ C_{\dar\uar
3}=C_{\dar\uar 4}=\emptyset\\
\end{array}$$
\end{ex}

\begin{rem} Again note that actions $\uar\dar$ and $\dar\uar$ are inverse,
that is if $K_{(i_s,j_s)}(\sigma)\ne\emptyset$ and
$F_{(i_s\rightleftarrows i_t)}(\sigma)\in C_{i_s\uar\dar}(\sigma)$
then $F_{(i_s\rightleftarrows i_t)}(F_{(i_s\rightleftarrows
i_t)}(\sigma)=\sigma\in C_{\dar\uar i_s}(F_{(i_s\rightleftarrows
i_t)}(\sigma));$ if $L_{(i_t,j_t)}(\sigma)\ne \emptyset$ and
$F_{(i_s\rightleftarrows i_t)}(\sigma)\in C_{\dar\uar i_s}(\sigma)$
then $F_{(i_s\rightleftarrows i_t)}(F_{(i_s\rightleftarrows
i_t)}(\sigma))=\sigma\in C_{i_s\uar\dar}(F_{(i_s\rightleftarrows
i_t)}(\sigma)).$
\end{rem}

\subsection{}\label{4.7}
By \cite[3.15]{Mel1} one has

\begin{thm}\label{desc}
 For $\sigma=(i_1,j_1)\ldots(i_k,j_k)\in\bS_n^2$ one has
$$\begin{array}{ll}
D(\sigma)=&\coprod\limits_{s=1}^k\sigma_{i_s\dar}\cup\coprod
\limits_{s=1}^k \sigma_{j_s\rar}\cup \coprod\limits_{s=2}^k
C_{i_s\looparrowright}(\sigma)\cup \coprod\limits_{s=1}^k
C_{i_s\uar\dar}(\sigma)\\
A(\sigma)=&\coprod\limits_{s=1}^k\sigma_{\uar i_s}\cup\coprod
\limits_{s=1}^k \sigma_{\lar j_s}\cup \coprod\limits_{s=2}^k
C_{\looparrowleft i_s}(\sigma)\cup \coprod\limits_{s=1}^k
C_{\dar\uar i_s}(\sigma)\\
\end{array}$$
\end{thm}

\begin{proof}
The equality for $D(\sigma)$ is exactly the contents of theorem 3.15
of \cite{Mel1}.

The equality for $A(\sigma)$ is a straightforward corollary of this
theorem. Indeed, $\sigma'\in A(\sigma)$ iff $\sigma\in D(\sigma').$
Thus, since our ``left" and "right" actions are inverse we get the
equality.
\end{proof}
\subsection{}\label{4.8}
For $\sigma,\sigma'\in\bS_n^2(k).$ it is obvious by the definition
of $D(\sigma)$ that
 $\codim_{\overline \B_\sigma}\B_{\sigma'}=1$ implies $\sigma'\in
 D(\sigma).$ Moreover, one has
\begin{thm}\label{thm3} Let $\sigma,\sigma'\in \bS_n^2(k)$ be  such that $\sigma'<\sigma.$
Then $\codim_{\overline\B_\sigma}\B_{\sigma'}=1$ if and only if
$\sigma'\in D(\sigma).$
\end{thm}
\begin{proof}
By the note at the beginning of the section it is obvious that if
$\sigma,\sigma'\in\bS_n^2(k)$ are such that
$\codim_{\overline\B_\sigma}\B_{\sigma'}=1$ then $\sigma'\in
D(\sigma).$

We have to show that for $\sigma'\in D(\sigma)$ one has
$\codim_{\overline\B_\sigma}\B_{\sigma'}=1.$ We do this by
straightforward computation for each type of elements in
$D(\sigma).$

 (i)\ \ \  Let
$\sigma'=\sigma_{i_s\dar}=(i_1,j_1)\ldots(m,j_s)\ldots(i_k,j_k).$
Recall $q_{(i_t,j_t)}(\sigma)$ defined by \ref{2.3}$(*)$. Note that
$q_{(i_t,j_t)}(\sigma_{i_s\dar})=q_{(i_t,j_t)}(\sigma)$ for any
$t\ne s.$ If $m=i_s-1$ then
$q_{(m,j_s)}(\sigma_{i_s\dar})=q_{(i_s,j_s)}(\sigma).$ If $m=i_s-p$
where $p>1$ then
$$
\begin{array}{ll}
q_{(i_s,j_s)}(\sigma)=&\#\{i_p<m\ :\ j_p<j_s\}+\#\{j_p\ :\
j_p<m\}+\#\{m<i_p<i_s\ :\ j_p<j_s\}\\
& +\#\{j_p\ :\ m<j_p<i_s\}\end{array}$$
and
$$q_{(m,j_s)}(\sigma_{i_s\dar})=\#\{i_p<m\ :\
j_p<j_s\}+\#\{j_p\ :\ j_p<m\} .$$
\\By Remark \ref{remark} we have
$\#\{m<i_p<i_s\ :\ j_p<j_s\}+\#\{j_p\ :\ m<j_p<i_s\}=i_s-m-1$ so
that
$q_{(m,j_s)}(\sigma_{i_s\dar})=q_{(i_s,j_s)}(\sigma)-(i_s-(m-1)).$
Thus,
$$\begin{array}{ll}
\dim\B_{\sigma'}&=\left\{\begin{array}{ll} kn-\sum\limits_{t=1,t\ne
s}^k(j_t-i_t)-(j_s-i_s+1)-\sum\limits_{t=2}^k q_t(\sigma) & {\rm
if}\
m=i_s-1\\
\begin{array}{l}
kn-\sum\limits_{t=1,t\ne
s}^k(j_t-i_t)-(j_s-m)\\
\quad\quad -\sum\limits_{t=2,t\ne
s}^k q_t(\sigma)-q_s(\sigma)+(i_s-(m-1))\\
\end{array}
 & {\rm otherwise}\\
\end{array}\right.\\
\quad&\quad \\
 &=\dim\B_{\sigma}-1\ .\\ \end{array}$$

(ii) Assume $\sigma_{j_s\rar}\ne\emptyset.$ To show
$\codim_{\overline\B_\sigma}\B_{\sigma_{j_s\rar}}=1$ we use the
symmetry about $y=n-x$ for the points of $\sigma.$ Indeed, let
$\hat\sigma=(n+1-j_1,n+1-i_1)\ldots(n+1-j_k,n+1-i_k)$ and put $(\hat
i_s,\hat j_s)=(n+1-j_s,n+1-i_s).$  Note that ${\hat\sigma}$ as a set
of points is obtained from $\sigma$ by symmetry about $y=n-x$.
Respectively $\hat\sigma_{\hat i_s\dar}$ is obtained from
$\sigma_{j_s\rar}$ by symmetry about $y=n-x.$ So by (i) we get
$\codim_{\overline\B_\sigma}\B_{\sigma_{j_s\rar}}=
\codim_{\overline\B_{\hat\sigma}}\B_{\hat\sigma_{\hat i_s\dar}}=1$

(iii)\ \ \ To show the claim for an element of type (II)-a) assume
that $C_{i_t\looparrowright}(\sigma)\ne\emptyset$  and let
$\sigma'=F_{(j_s\rightleftarrows i_t)}(\sigma)\in
C_{i_t\looparrowright}(\sigma)$ We write the new pairs at the place
 of the old ones so that
 $$\sigma'=(i_1,j_1)\ldots(i_s,i_t)\ldots(j_s,j_t)\ldots.$$
 If $i_t=j_s+1$ then $\sigma'$ is written in the canonical form. Otherwise pair
 $(j_s,j_t)$ can
 be not on the right place. But in any case note that for any $p\ne s,t$ the
 contribution to $q_p(\sigma)$ of all pairs but pairs $s,t$ is equal to their
 contribution to $q_p(\sigma').$ Thus,
 \begin{itemize}
 \item[(a)] For any $p\ :\ p\leq s$  one has
 $q_p(\sigma')=q_p(\sigma).$ Indeed, this is obviously true for
 $p<s$. As well, $q_s(\sigma')=q_s(\sigma)$ since by
 Remark \ref{rem1}(i) we have that for any $p<s$ either $j_p<j_s$
 or $j_p>i_t.$

 \item[(b)] For any $p\ :\ p>t$ one has
 $q_p(\sigma')=q_p(\sigma)$. Indeed, if $i_p>i_t$ then the pair $(i_s,j_s)$
 adds 2 to $q_p(\sigma)$ and the pair $(i_s,i_t)$ adds 2, the pair
 $( i_t,j_t)$ adds the same number to $q_p(\sigma)$ as the pair
 $(j_s,j_t)$ to $q_p(\sigma').$
 \item[(c)] For any $p\ :\ s<p<t$ and $j_p<j_s$ one has $q_p(\sigma)=q_p(\sigma').$
 Indeed, note first that $i_p<j_s$ so that the contribution of
 pairs $(i_s,j_s)$ and $(i_t,j_t)$ to $q_p(\sigma)$ is 0 as well as
 the contribution of pairs $(i_s,i_t)$ and $(j_s,j_t)$ to
 $q_p(\sigma').$
 \item[(d)] For any $p\ :\ s<p<t$ and $j_p>j_t$ one has
 $q_p(\sigma)=q_p(\sigma').$ Indeed, by Remark \ref{rem1}(ii) these means that
   $i_p<j_s$ so that  the contribution of $(i_s,j_s)$
 to $q_p(\sigma)$ is 1 and the contribution  of $(i_s,i_t)$ to
 $q_p(\sigma')$ is 1; the contribution of $(i_t,j_t)$ to
 $q_p(\sigma)$ is 0 and the contribution of $(j_s,j_t)$ to
 $q_p(\sigma')$ is 0.
 \end{itemize}

Summarizing, (a)-(d) provide us that for  $p\ne t$ such that
$\{i_p,j_p\}\cap\{a\}_{a=j_s+1}^{i_t-1}=\emptyset$ one has
$q_p(\sigma)=q_p(\sigma').$

 Moreover by Remark \ref{rem1}(iii) if $p$ is
such that $\{i_p,j_p\}\cap\{a\}_{a=j_s+1}^{i_t-1}\ne\emptyset$ then
$i_s<i_p<j_p<j_t.$ Let us consider $q_p$ for such $p.$
 If $i_t=j_s+1$ this set is empty.
 Assume $i_t=j_s+b$ where $b>1$

 \begin{itemize}
 \item[(1)] Assume  $i_p\in \{a\}_{a=j_s+1}^{i_t-1}$ for some
 $p\ :\ s<p<t$. Then
 \begin{itemize}
 \item[(a)] if $j_s<j_p<i_t$ then $q_p(\sigma')=q_p(\sigma)-2;$
 \item[(b)] if $j_p>i_t$ then $q_p(\sigma')=q_p(\sigma)-1;$
\end{itemize}
\item[(2)] Assume $j_p\in \{a\}_{a=j_s+1}^{i_t-1}$ for some $p\ :\ s<p<t.$
The case $j_s<i_p<i_t$ is already considered in (i(a)). So we have
to consider only the case when $i_p<j_s.$ In that case
$q_p(\sigma')=q_p(\sigma)-1.$
\end{itemize}
Summarizing we get
$$\sum\limits_{p=1,p\ne t}^{k}q_p(\sigma')=\sum\limits_{p=1,p\ne t}^{k}q_p(\sigma)-(i_t-j_s-1)$$
Finally,
$q_t(\sigma')=q_t(\sigma)-(i_t-j_s-1)-1=q_t(\sigma)-(i_t-j_s).$
Indeed for any $p\ :\  p<s$ one has by Remark \ref{rem1}(i) that
either $j_p<j_s$ or $j_p>i_t$   so that its contribution to
$q_t(\sigma)$ is either 2 or 0 respectively and correspondingly its
contribution to $q_t(\sigma)$ is either 2 or 0 respectively.
Further, the contribution of $(i_s,j_s)$ to $q_t(\sigma)$ is 2 and
contribution of $(i_s,i_t)$ to $q_t(\sigma')$ is 1. For any $p\ :\
s<p<t$ one of three following situations is possible:
\begin{itemize}
\item[(a)] $\{i_p,j_p\}\cap\{a\}_{a=j_s+1}^{i_t-1}=\emptyset.$ In this
case the contribution of $(i_p,j_p)$ to $q_t(\sigma)$ is equal to
the contribution of $q_t(\sigma');$
\item[(b)] $\{i_p,j_p\}\cap\{a\}_{a=j_s+1}^{i_t-1}\ne\emptyset$ and
$\{i_p,j_p\}\not\subset\{a\}_{a=j_s+1}^{i_t-1}.$ In this case the
contribution of $(i_p,j_p)$ to $q_t(\sigma)$ is 1 more than its
contribution to $q_t(\sigma').$
\item[(c)] $\{i_p,j_p\}\subset\{a\}_{a=j_s+1}^{i_t-1}.$ Then the
contribution of $(i_p,j_p)$ to $q_t(\sigma)$ is 2 and its
contribution to $q_t(\sigma')$ is 0.
\end{itemize}
Summarizing, any  $a\ :\ j_s<a<i_t$ adds 1 to $q_t(\sigma)$
comparing with $q_t(\sigma')$.

Now we can compute $\dim(\B_{\sigma'}):$
$$\begin{array}{ll}
   \dim(\B_{\sigma'})&=\left\{
             \begin{array}{ll}
              \begin{array}{l}
                   kn-\sum\limits_{p=1,p\ne s,t}^k
                            (j_p-i_p)-(j_s+1-i_s)-(j_t-i_t+1)\\
                   -\sum\limits_{p=1,p\ne t}^k q_p(\sigma)-(q_t(\sigma)-1)\\
              \end{array}
         & {\rm if}\ i_t=j_s+1\\
              \begin{array}{l}
               kn-\sum\limits_{p=1,p\ne s,t}^k(j_p-i_p)-(i_t-i_s)-(j_t-j_s)\\
               -\sum\limits_{p=1,p\ne t}^{k}q_p(\sigma)+(i_t-j_s-1)
               -(q_t(\sigma)-(i_t-j_s))\\
               \end{array}
               &{\rm otherwise}\\
\end{array}\right.\\
\quad&\quad\\
&=\dim(\B_{\sigma})-1\ .\\
\end{array}$$

(iv) Finally, assume $C_{i_s\uar\dar}(\sigma)\ne \emptyset$ and  let
$\sigma'=F_{(i_s \rightleftarrows i_t)}(\sigma)\in
C_{i_s\uar\dar}(\sigma).$ Again, we write the new pairs at the place
 of the old ones so that $\sigma'=(i_1,j_1)\ldots(i_s,j_t)\ldots(i_t,j_s)\ldots.$
 Note that $\sigma'$ is written in canonical form.
 Exactly as in the previous
case $q_p(\sigma')=q_p(\sigma)$ for any $p\ :\ p<s$ or $p>t.$ As
well by $(*)$ from \ref{4.64} one has $ q_p(\sigma')=q_p(\sigma)$
for any $p\ :\ s<p<t.$ Now let us compute
$q_s(\sigma')+q_t(\sigma').$ First, note that for any $p<s$ such
that $j_s<j_p<j_t$ we have to subtract 1 from $q_t(\sigma')$ and to
add 1 to $q_s(\sigma').$ We have to add 1 to $q_t(\sigma')$ since
$j_t<j_s.$ As well for any $p\ :\ s<p<t$ one has by $(*)$ from
\ref{4.64} that $j_p<j_t$ if and only if $j_p<j_s.$ Summarizing, we
get $q_s(\sigma')+q_t(\sigma')=q_s(\sigma)+q_t(\sigma)+1.$ Thus,
\begin{eqnarray*}\dim(\B_\sigma')&=& kn-\sum\limits_{p=1,p\ne
s,t}^k(j_p-i_p)-(j_t-i_s)-(j_s-i_t)\\
&&-\sum\limits_{p=1,p\ne s,t}^k
q_p(\sigma)-(q_s(\sigma)+q_t(\sigma)+1)\\
&=&\dim(\B_\sigma)-1.
\end{eqnarray*}
\end{proof}

Note that our proof is a somewhat more technical than
 the proofs
of lemmas \cite[3.11, 3.13, 3.14]{Mel1}, However it  shows a little
bit more than these lemmas. It shows that the set denoted by
$D_1(\sigma)$ in \cite{Mel1} is equidimensional of codimension 1 in
$\overline\B_\sigma.$

\subsection{}\label{4.9}
Let us discuss briefly $(\bS_n^2,\prec)$ as a poset. Note that $Id$
is the minimal  element of $\bS_n^2$ and it is in order with all the
elements of $\bS_n^2.$ By Remark \ref{rem2} (i)
$\sigma_o(k):=(1,n-k+1)\ldots(k,n)$ is the minimal element in
$\bS_n^2(k)$, it is easy to see that for $j\geq k$ we have
$\sigma_o(j)\succeq\sigma_o(k)$; in particular we deduce that
$\sigma_o(k)$ is the minimal element in $\bigsqcup_{j\geq
k}\bS_n^2(j)$. We can define a generalized descendant of $\sigma$ in
$\bS_n^2$ generalizing the notion of descendant in
$\bS_n^2(L(\sigma)).$ Namely put $\sigma'$ to be a generalized
descendant of $\sigma$ if $\sigma'\prec \sigma$ and
$\sigma'\preceq\phi\preceq\sigma$ implies either $\phi=\sigma$ or
$\phi=\sigma'.$ The set $C(\sigma)$ of generalized descendants is
what is called usually the cover of $\sigma$ for the order
$\preceq.$

Note that $C(\sigma)\cap \bS_n^2(L(\sigma))=D(\sigma)$ so that it is
 described in Theorem \ref{desc}. For $\sigma'\prec \sigma$ such that
and $L(\sigma')<L(\sigma)$ by Proposition 3.9 of \cite{Mel1} there
exists $\sigma"$ obtained from $\sigma$ by deleting one point
$P:(i_s,j_s)\in\sigma$ such that $\sigma'\preceq\sigma''$. Moreover
$\sigma"$ is maximal possible if for any $(i_t,j_t)\in\sigma$ one
has that $i_t<i_s$ implies $j_t<j_s$ (one also can see this easily
from Lemma \ref{lem1}.)

 Let us first define
\begin{defn}
\begin{enumerate}
\item[(i)] Let $\sigma\succ\sigma'.$
Call a {\bf chain} between $\sigma$ and $\sigma'$ a sequence of
generalized descendants that is
$$\sigma=\sigma_m\succ
\sigma_{m-1}\succ\sigma_{m-2}\succ\cdots\succ\sigma_0=\sigma'$$
 where $\sigma_{i-1}$ is a generalized descendant of $\sigma_i$ for
$i\ :\ 1\leq i\leq m$

 The number $m$ is called the {\bf length} of the chain.
\item[(ii)] For $\sigma\in\bS_n^2(j)$ where $j\geq k$ put
${\rm depth}_k(\sigma)$ to be the length of chain between $\sigma$
and $\sigma_o(k).$
\end{enumerate}
\end{defn}

By Theorem \ref{thm3} for $\sigma, \sigma' \in \bS_n^2(k)$ such that
$\sigma'\prec\sigma$ one has that the
 length of the chain between $\sigma$ and $\sigma'$ does not depend on the
 choice of the chain and moreover it equals to $\codim_{\ov\B_\sigma}(\B_{\sigma'}).$

Now we can state a more general result:

\begin{thm}\label{lem2}
Let $\sigma'$ be a generalized descendant of $\sigma$. Then
$\codim_{\overline\B_\sigma}\B_{\sigma'}=1$. Thus,
\begin{enumerate}
\item[(i)] For  $\sigma\succ\sigma'$
any chain between $\sigma'$ and $\sigma$ is of length equal to
$\codim_{\ov\B_\sigma}(\B_{\sigma'}).$
\item[(ii)] In particular, for any $\sigma\in \bS_n^2(l)$ where
$l\geq k$ ${\rm depth}_k(\sigma)$ is well defined and  ${\rm
depth}_k(\sigma)=\codim_{\ov\B_\sigma}\B_{\sigma_o(k)}.$
\item[(iii)] In particular, ${\rm depth}_0(\sigma)=\dim \B_\sigma$.
\end{enumerate}
\end{thm}

\begin{proof}
The only case that we have to check is the case where $\sigma'$ a
generalized descendant of $\sigma$ obtained by deleting a point in
its support.

For any point $P=(x,y)$, denote
$\pspolygon[fillstyle=solid,fillcolor=gray](0,0.35)(0.5,0.35)(0.5,-0.15)$\hspace{0.6cm}$(P)$
and,
$\pspolygon[fillstyle=solid,fillcolor=lightgray](0,0.35)(0.5,0.35)(0.5,-0.15)(0,-0.15)$\hspace{0.6cm}$(P)$
the areas drawn in Fig. \ref{Fig6}.

\begin{figure}[htbp]
  \centering
\psset{unit=0.5cm}
\begin{tabular}{cc}
  \begin{pspicture}(-1,-3)(6,0)

\psline(0,0)(4,-4) \put(0.5,-2){$\Delta$}

\psline[linestyle=dotted,dotsep=1pt](5,0)(5,-4)
\psline[linestyle=dotted,dotsep=1pt](5.5,-2)(2,-2)(2,0)

\pspolygon[fillstyle=solid,fillcolor=lightgray](5,-2)(2,-2)(2,0)(5,0)
\put(4.8,-2.2){$\bullet$}\put(5.2,-1.8){$P$}

\pspolygon[fillstyle=solid,fillcolor=gray](2,-2)(2,0)(0,0)

\psline[linestyle=dotted,dotsep=1pt](5.5,-1)(1,-1)(1,0)

\end{pspicture}
\end{tabular}\caption{} \label{Fig6}
\end{figure}

\noindent Any point $Q=(x',y')$ strictly contained in
$\pspolygon[fillstyle=solid,fillcolor=gray](0,0.35)(0.5,0.35)(0.5,-0.15)$\hspace{0.6cm}$(P)$
(resp. in
$\pspolygon[fillstyle=solid,fillcolor=lightgray](0,0.35)(0.5,0.35)(0.5,-0.15)(0,-0.15)$\hspace{0.6cm}$(P)$)
verifies $x'<x$ and $y'<x$ (resp $x'<x$ and $x<y'$).

 (a) Let $P_o=(i_o,j_o)$ be the point deleted from
$\sigma$ to get $\sigma'$. If $\{r\}_{r=1}^{i_o-1}\cap
\langle\sigma\rangle^c\neq\emptyset$, then $\sigma_{\uparrow i_o}$
exists and $\sigma'\prec\sigma_{\uparrow i_o}\prec\sigma$ which is
impossible by choice of $\sigma'$. Therefore any integer
$x\in\{r\}_{r=1}^{i_o-1}$ is either a first entry or a second entry
of a point of $N_\sigma$. Notice that if $x$ is a second entry,
since there is no point of $N_\sigma$ in the North-East of $P$, the
corresponding point in $N_\sigma$ has to be in
$\pspolygon[fillstyle=solid,fillcolor=gray](0,0.35)(0.5,0.35)(0.5,-0.15)$\hspace{0.6cm}$(P_o)$,
see Fig. \ref{Fig6}.

\noindent Denote $Q_1,\dots,Q_l)$ (resp. $R_1,\cdots,R_m$) the
points of $N_\sigma$ inside of
$\pspolygon[fillstyle=solid,fillcolor=gray](0,0.35)(0.5,0.35)(0.5,-0.15)$\hspace{0.6cm}$(P_o)$
(resp.
$\pspolygon[fillstyle=solid,fillcolor=lightgray](0,0.35)(0.5,0.35)(0.5,-0.15)(0,-0.15)$\hspace{0.6cm}$(P_o)$.
In particular, we get
\begin{eqnarray}\label{eq3}
2l+m=i_o-1\ {\rm and}\ q_{P_o}(\sigma)= i_o-1.
\end{eqnarray}
Moreover for those points which are in the North-West of $P_o$, we
have
\begin{eqnarray}\label{eq4}
q_{Q_i}(\sigma')=q_{Q_i}(\sigma)\ {\rm and}\
q_{R_i}(\sigma')=q_{R_i}(\sigma').
\end{eqnarray}

(b) In the same way, we have $\{r\}_{r=j_o+1}^{n}\cap
\langle\sigma\rangle^c=\emptyset$, and any integer
$y\in\{r\}_{r=j_o+1}^{n}$ is either a first entry or a second entry
of a point of $N_\sigma$; if $y$ is a first entry of a point $Q$ of
$N_\sigma$ then we have $P_o\in
\pspolygon[fillstyle=solid,fillcolor=gray](0,0.35)(0.5,0.35)(0.5,-0.15)$\hspace{0.6cm}$(Q)$.
Denote $S_1,\cdots,S_u$ (resp. $T_1,\cdots,T_v$) the points of
$N_\sigma$ such that $P_o$ lies inside
$\pspolygon[fillstyle=solid,fillcolor=gray](0,0.35)(0.5,0.35)(0.5,-0.15)$\hspace{0.6cm}$(S_i)$
(resp.
$\pspolygon[fillstyle=solid,fillcolor=lightgray](0,0.35)(0.5,0.35)(0.5,-0.15)(0,-0.15)$\hspace{0.6cm}$(T_i)$,
then we have
\begin{eqnarray}\label{eq5}
2u+v=n-j_o.
\end{eqnarray}
 And we have the relations
\begin{eqnarray}\label{eq6}
q_{S_i}(\sigma')=q_{S_i}(\sigma)-2\ {\rm and}\
q_{T_i}(\sigma')=q_{T_i}(\sigma')-1.
\end{eqnarray}

(c) Finally any point $Q=(i_s,j_s)$ in the South-West of $P$
verifies \begin{eqnarray}\label{eq7}q_{Q}(\sigma')=q_{Q}(\sigma).
\end{eqnarray}

\noindent The formula in Theorem \ref{dim} does not depend on the
ordering of the cycles in $\sigma$, but depends just on the
distribution of points of $N_\sigma$, therefore let us put
$P_s=(i_s,j_s)$ for $1\leq s\leq k$ [of course $P_o$ still
represents one of them], and the notation ``$P_s\ {\rm in}\ SE$"
(resp. ``$P_s\ {\rm not}\ SE$") will mean that $P_s$ is in the
South-East (resp. not in the South-East) of the point $P_o$.

then we have
$$
\dim \B_\sigma= kn-\sum_{P_s\in N_\sigma}(j_s-i_s)-\sum_{P_s\neq
P_1}q_{P_s}(\sigma)$$
$$=(k-1)n-\sum_{P_s\in N_\sigma-\{P_o\}}(j_s-i_s)+(n-j_o+i_o)-\sum_{\substack {P_s\neq P_1\\
P_s\ {\rm not}\ SE}}q_{P_s}(\sigma)-\sum_{P_s\ {\rm in}\
SE}q_{P_s}(\sigma)-q_{P_o}(\sigma). $$

\noindent From (\ref{eq3}), (\ref{eq4}) and (\ref{eq7}) we have
\begin{eqnarray}\label{eq8}
q_{P_o}(\sigma)= i_o-1\ {\rm and}\ \sum\limits_{\substack {P_s\neq P_1\\
P_s\ {\rm not}\ SE}}q_{P_s}(\sigma)=\sum\limits_{\substack {P_s\neq P_1\\
P_s\ {\rm not}\ SE}}q_{P_s}(\sigma').
\end{eqnarray}

\noindent On the other hand by (\ref{eq5}) and (\ref{eq6}) we have
\begin{eqnarray}\label{eq9}
\sum\limits_{P_s\ {\rm in}\ SE}q_{P_s}(\sigma)=\sum\limits_{P_s\
{\rm in}\ SE}q_{P_s}(\sigma')+(n-j_o).
\end{eqnarray}

\noindent Now combining (\ref{eq8}) and (\ref{eq9}) we get

$$\dim \B_{\sigma}=\underbrace{(k-1)n-\sum_{P_s\in N_\sigma-\{P_o\}}(j_s-i_s)-\sum
\limits_{P_s\neq P_o,P_1}q_{P_s}(\sigma')}_{\dim \B_{\sigma'}}+1.$$

\noindent and the proof is done.

\end{proof}

\subsection{}\label{4.9}
Now let us return to the intersections of the closures of
codimension 1.

 Note that for any $\sigma',\sigma\prpr\in A(\sigma)$
one has
$\B_\sigma\in\overline\B_{\sigma'}\cap\overline\B_{\sigma\prpr}$
Thus by Theorem  \ref{thm3}
$\codim_{\overline\B_{\sigma'}}(\overline\B_{\sigma'}\cap\B_{\sigma\prpr})=1.$
Moreover by this theorem, if $\sigma,\sigma'\in \bS_n^2(k)$ are such
that $\dim\B_\sigma=\dim\B_{\sigma'}$ then
$\codim_{\overline\B_\sigma}(\overline\B_\sigma\cap\overline\B_{\sigma'}\cap\Or_{(n-k,k)})=1$
iff there exists $\sigma\prpr\in \bS_n^2(k)$ such that
$\sigma,\sigma\pr\in A(\sigma\prpr).$ In general this does not imply
that the intersection $\overline\B_\sigma\cap\overline\B_{\sigma'}$
is irreducible.

The first example of reducible intersection of codimension 1 occurs
in $n=5.$ Let $\sigma\prpr=(1,5)(2,4).$ Let $\sigma=(1,5)(3,4)$ and
$\sigma'=(2,4)(3,5).$ By Theorem \ref{thm3} $\sigma,\sigma'\in
A(\sigma\prpr)$ thus their intersection is of codimension 1. Let us
compute $R_{\sigma,\sigma'}.$

$$
R_{\sigma}=\left( \tiny{\begin{array}{ccccc}
0&0&0&1&2\\
0&0&0&1&1\\
0&0&0&1&1\\
0&0&0&0&0\\
0&0&0&0&0\\
\end{array}}
\right),\ R_{\sigma'}=\left( \tiny{\begin{array}{ccccc}
0&0&0&1&2\\
0&0&0&1&2\\
0&0&0&0&1\\
0&0&0&0&0\\
0&0&0&0&0\\
\end{array}}
\right) \ \Longrightarrow\
 R_{\sigma,\sigma'}=\left(\tiny{\begin{array}{ccccc}
0&0&0&1&2\\
0&0&0&1&1\\
0&0&0&0&1\\
0&0&0&0&0\\
0&0&0&0&0\\
\end{array}}\right)$$
Note that $R_{\sigma,\sigma'}\ne R_{\sigma''}$ so the intersecion is
reducible. Indeed, one can easily see that
$$\overline\B_\sigma\cap\overline\B_{\sigma'}=\overline\B_{(1,4)(3,5)}\cup
\overline\B_{(1,5)(2,4)}$$
 Note that $\dim\B_{(1,4)(3,5)}=\dim\B_{(1,5)(2,4)}=4$ so that the
 intersection is equidimensional.

\section{Intersections of codimension 1 of orbital varieties of
nilpotent order 2}
\subsection{}\label{5.1}
Now we apply the machinery developed in the previous section to
orbital varieties of nilpotent order 2. Our first aim is to show
that for $T,S\in \bT_{(n-k,k)}$ if $\codim_{\V_T}(\V_T\cap\V_S)=1$
then this intersection is irreducible. Since
$\codim_{\V_T}(\V_T\cap\V_S)=1$ if and only if $D(\sigma_T)\cap
D(\sigma_S)\ne\emptyset$ it is enough to consider $A(\sigma')$ for
any $\sigma'\in D(\sigma_T).$ We show

\begin{thm}\label{thm}
For any $k\leq \frac{n}{2}$ for any $T\in\bT_{(n-k,k)}$ for any
$\sigma'\in D(\sigma_T)$ one has $A(\sigma')$ contains only 2
elements namely $\{\sigma_T,\sigma_S\}$ where $S\in\bT_{(n-k,k)}$
and
$\ov\B_T\cap\ov\B_S\cap\Or_{(n-k,k)}=\ov\B_{\sigma'}\cap\Or_{(n-k,k)}$
so that this intersection is irreducible.
\end{thm}

\begin{proof}
We show this by considering all possible types of $\sigma'\in
D(\sigma_T).$

Put $\sigma=\sigma_T=(i_1,j_1)\ldots(i_k,j_k).$ By remark (ii) from
\ref{4.5} $\dim(\B_\sigma)=k(n-k)$ and is maximal possible. By
Theorem \ref{thm3} for any $\sigma'\in D(\sigma)$
$\dim(\B_{\sigma'})=k(n-k)-1.$
\begin{itemize}
\item[(i)] Assume $\sigma'=\sigma_{i_s\dar}\ne\emptyset$ and let
$\sigma'=F_{(i_s\rar m)}(\sigma).$ Let $\sigma\prpr=F_{(j_s\rar
i_s)}(\sigma')$ that is obtained from $\sigma$ by changing
$(i_s,j_s)$ to $(m,i_s).$ Since $m\not \in <\sigma>$ we get that
$\sigma\prpr$ is correctly defined. Moreover $\sigma\prpr>\sigma'$
so that $\dim\B_{\sigma\prpr}>\dim(\B_{\sigma'})$, thus,
$\dim(\B_{\sigma\prpr})=k(n-k)$ and $\sigma\prpr\in A(\sigma').$ It
is left to show that $R_{\sigma,\sigma\prpr}=R_{\sigma'}$ and that
$A(\sigma')=\{\sigma,\sigma\prpr\}.$

Let us first compute $R_{\sigma,\sigma\prpr}.$ By \ref{4.4} for any
$j<j_s$
$(R_\sigma)_{i,j}=(R_{\sigma_{(i_s,j_s)}^-})_{i,j}=(R_{\sigma'})_{i,j}$
and for any $i>m$
$(R_{\sigma\prpr})_{i,j}=(R_{\sigma_{(i_s,j_s)}^-})_{i,j}=(R_{\sigma'})_{i,j}.$
If $i\leq m$ and $j\geq j_s$ then
$(R_{\sigma'})_{i,j}=(R_\sigma)_{i,j}=(R_{\sigma\prpr})_{i,j}$ so
that $R_{\sigma,\sigma\prpr}=R_{\sigma'}.$

Assume that $w\ne\sigma,\sigma\prpr\in A(\sigma').$ In four items
below we show that $w$ cannot be one of four types of ancestors we
define, so that $A(\sigma')=\{\sigma,\sigma\prpr\}.$
\begin{itemize}
\item[(a)] Assume  $w=\sigma'_{\uar i_t}=F_{(i_t\rar a)}(\sigma').$
Then $i_t\ne m$ (otherwise $w=\sigma$).  If $a\ne i_s$ then
$F_{(i_t\rar a)}(\sigma)$ is well defined and $F_{(i_t\rar
a)}(\sigma)>\sigma$ which contradicts maximality of $\sigma.$ If
$a=i_s$ and $i_t>m$ then $i_s<j_t<j_s$ by definition of $\sigma'$ so
that $F_{(i_t\rightleftarrows i_s)}(\sigma)$ is well defined and
$F_{(i_t\rightleftarrows i_s)}(\sigma)>\sigma$ which contradicts
maximality of $\sigma.$ If $a=i_s$ and $i_t<m$ then $F_{(i_t\rar
m)}(\sigma)$ is well defined and $F_{(i_t\rar m)}(\sigma)>\sigma$
which contradicts maximality of $\sigma.$ Thus, $w\ne \sigma'_{\uar
i_t}.$
\item[(b)] Exactly in the same way comparing $w$ with $\sigma\prpr$
we get by maximality of $\sigma\prpr$ that $w\ne \sigma'_{\lar
j_t}.$
\item[(c)] Assume $w=F_{(i_p\leftrightarrows i_q)}(\sigma')$ where
$i_p<i_q$ so that $i_q<j_p<j_q.$
If $i_p, i_q\ne m$ then $F_{(i_p\leftrightarrows i_q)}(\sigma)$ is
well defined and $F_{(i_p\leftrightarrows i_q)}(\sigma)>\sigma$
which contradicts to maximality of $\sigma.$ If $i_p=m$ then by
definition of $\sigma'$ we get that $i_q>i_s$ thus
$F_{(i_s\rightleftarrows i_q)}(\sigma)$ is well defined and
$F_{(i_s\rightleftarrows i_q)}(\sigma)>\sigma$ which contradicts to
maximality of $\sigma.$ Finally, if $i_q=m$ then $F_{(i_p\rar
m)}(\sigma)$ is well defined and $F_{(i_p\rar m)}(\sigma)>\sigma$
which contradicts to maximality of $\sigma.$ Thus, $w\ne
F_{(i_p\rightleftarrows i_q)}(\sigma').$
\item[(d)] Assume $w=F_{(i_q\rightleftarrows j_p)}(\sigma')$ where
$i_p<i_q<j_p<j_q.$ Again if
$i_p,i_q\ne m$ then $F_{(i_q\rightleftarrows j_p)}(\sigma)$ is well
defined and $F_{(i_q\rightleftarrows j_p)}(\sigma)>\sigma$ which
contradicts to maximality of $\sigma.$ If $i_q=m$ then $F_{(j_p\rar
m)}(\sigma)$ is well defined and $F_{(j_p\rar m)}(\sigma)>\sigma$
which contradicts to maximality of $\sigma.$ Finally, if $i_p=m$
then $i_q>i_s$ since by definition of $\sigma'$ for all $i_r\ :\
m<i_r<i_s$ one has $j_r<j_s.$ Thus we get $i_s<i_q<j_s<j_q$ so that
$F_{(i_q\rightleftarrows j_s)}(\sigma)$ is well defined and
$F_{(i_q\rightleftarrows j_s)}(\sigma)>\sigma$ which contradicts to
maximality of $\sigma.$ Thus, $w\ne F_{(i_q\rightleftarrows
j_p)}(\sigma').$
\end{itemize}

\item[(ii)] Assume $\sigma'=\sigma_{ j_s\rar}\ne\emptyset$ and let
$\sigma'=F_{(j_s\rar m)}(\sigma).$ Let $\sigma\prpr=F_{(i_s\rar
j_s)}(\sigma')$ that is obtained from $\sigma$ by changing
$(i_s,j_s)$ to $(j_s,m).$ One has exactly as in (i) that
$\sigma\prpr$ is maximal and that $\sigma'=\sigma\prpr_{j_s\dar}.$
Thus by (i) $A(\sigma')=\{\sigma,\sigma\prpr\}$ and
$R_{\sigma,\sigma\prpr}=R_{\sigma'}.$
\item[(iii)] Assume $C_{i_s\uar\dar}(\sigma)\ne\emptyset$ and let
$F_{(i_s \rightleftarrows i_t)}(\sigma)\in C_{i_s\uar\dar}(\sigma).$
Recall that this implies $i_s<i_t<j_t<j_s.$ Let
$\sigma\prpr=F_{(i_t\rightleftarrows j_t)}(\sigma')$ that is
obtained from $\sigma$ by changing pairs $(i_s,j_s)(i_t,j_t)$ to
pairs $(i_s,i_t)(j_t,j_s)$. Again $\sigma\prpr>\sigma'$ so that
$\dim\B_{\sigma\prpr}>\dim(\B_{\sigma'})$, thus,
$\dim(\B_{\sigma\prpr})=k(n-k)$ and $\sigma\prpr\in A(\sigma').$

Let us compute $R_{\sigma,\sigma\prpr}.$ Note that for $i\leq i_s,\
j\geq j_t$ and for $i\leq i_t,\ j\geq j_s$ one has
$(R_\sigma)_{i,j}=(R_{\sigma'})_{i,j}=(R_{\sigma\prpr})_{i,j}.$
Further, by lemma \ref{4.4} for $j<j_t$ or $i>i_t$
$(R_{\sigma})_{i,j}=(R_{\sigma_{(i_t,j_t)}^-})_{i,j}=(R_{\sigma'})_{i,j}.$
Finally for $i>i_s,\ j<j_s$
$(R_{\sigma\prpr})_{i,j}=(R_{\sigma_{i_t,j_t}^-})_{i,j}=(R_{\sigma'})_{i,j}.$
All this together provides $R_{\sigma,\sigma\prpr}=R_{\sigma'}.$

Assume that $w\ne\sigma,\sigma\prpr\in A(\sigma').$ In four items
below we show that $w$ cannot be one of four types of ancestors we
define, so that $A(\sigma')=\{\sigma,\sigma\prpr\}.$
\begin{itemize}
\item[(a)] Assume $w=\sigma'_{\uar i_p}=F_{(i_p\rar a)}(\sigma').$ Note that
for any $i_p\ne i_t$ $F_{(i_p\rar a)}(\sigma)$ is well defined and
$F_{(i_p\rar a)}(\sigma)>\sigma$ which contradicts maximality of
$\sigma.$ If $i_p=i_t$ then $F_{(i_s\rar a)}(\sigma)$ is well
defined and $F_{(i_s\rar a)}(\sigma)>\sigma$ which contradicts
maximality of $\sigma.$ Thus $w\ne\sigma'_{\uar i_p}.$
\item[(b)] By the symmetry around the anti-diagonal (i) implies $w\ne\sigma'_{\lar
j_p}.$
\item[(c)] Assume $w=F_{(i_p\rightleftarrows i_q)}(\sigma')$ where
$i_p<i_q$ so that $i_q<j_p<j_q.$ If
$\{i_p,i_q\}\cap\{i_s,i_t\}=\emptyset$ then $F_{(i_p\rightleftarrows
i_q)}(\sigma')$ is well defined and $F_{(i_p\rightleftarrows
i_q)}(\sigma')>\sigma$ which contradicts maximality of $\sigma.$ If
$i_p=i_s$ and $i_q<i_t$ then by definition of $\sigma'$ we get that
$j_q>j_s$ so that $F_{(i_s\rightleftarrows i_q)}(\sigma)$ is well
defined and $F_{(i_s\rightleftarrows i_q)}(\sigma)>\sigma$ which
contradicts maximality of $\sigma.$ If $i_p=i_s$ and $i_q>i_t$ then
$F_{(i_t\rightleftarrows i_q)}(\sigma)$ is well defined and
$F_{(i_t\rightleftarrows i_q)}(\sigma)>\sigma$ which contradicts
maximality of $\sigma.$ If $i_p=i_t$ then $j_q>j_s$ so that
$F_{(i_s\rightleftarrows i_q)}(\sigma)$ is well defined and
$F_{(i_s\rightleftarrows i_q)}(\sigma)>\sigma$ which contradicts
maximality of $\sigma.$ By symmetry around the anti-diagonal we get
that $i_q\ne i_s,i_t.$ Thus $w\ne F_{(i_p\rightleftarrows
i_q)}(\sigma').$
\item[(d)] Finally, assume $w=F_{(i_q\rightleftarrows j_p)}(\sigma')$
where $i_p<i_q<j_p<j_q.$
Again if $\{q,p\}\cap\{s,t\}=\emptyset$ then
$F_{(i_q\rightleftarrows j_p)}(\sigma)$ is well defined and
$F_{(i_q\rightleftarrows j_p)}(\sigma)>\sigma$ which contradicts
maximality of $\sigma.$ If $i_p=i_s$ and $i_q<i_t$ then by
definition of $\sigma'$ $j_q>j_s$ so that $F_{(i_q\rightleftarrows
j_s)}(\sigma)$ is well defined and $F_{(i_q\rightleftarrows
j_s)}(\sigma)>\sigma$ which contradicts maximality of $\sigma.$ If
$i_p=i_s$ and $i_q>i_t$ then $F_{(i_q\rightleftarrows j_t)}(\sigma)$
is well defined and $F_{(i_q\rightleftarrows j_t)}(\sigma)>\sigma$
which contradicts maximality of $\sigma.$ If $i_q=i_s$ then
$F_{(i_s\rightleftarrows j_p)}(\sigma)$ is well defined and
$F_{(i_s\rightleftarrows j_p)}(\sigma)>\sigma$ which contradicts
maximality of $\sigma.$ By the symmetry around the anti-diagonal we
get that $i_p,i_q\ne i_t.$
\end{itemize}
\item[(vi)] Assume $C_{i_t\looparrowright}(\sigma)\ne\emptyset$ and let
$\sigma'=F_{(j_s\rightleftarrows i_t)}(\sigma)\in
C_{i_t\looparrowright}(\sigma).$ Let
$\sigma\prpr=F_{(i_s\rightleftarrows j_s)}(\sigma')$ that is
obtained from $\sigma$ by changing $(i_s,j_s)(i_t,j_t)$ to
$(i_s,j_t)(j_s,i_t)$. One has exactly as in (iii) that
$\sigma\prpr\in C_{\dar\uar i_s}(\sigma')$ and $\sigma\prpr$ is
maximal. Thus by (iii) $A(\sigma')=\{\sigma,\sigma\prpr\}$ and
$R_{\sigma,\sigma\prpr}=R_{\sigma'}.$
 \end{itemize}
\end{proof}
\subsection{}\label{5.2}
Let us translate the results of the previous subsection into the
combinatorics of $S,T\in\bT_{(n-k,k)}$ such that
$\codim_{\V_T}(\V_T\cap\V_S)=1.$ This combinatorics works through
$\sigma_T.$

We need some notation. Let $k\leq\frac{n}{2}.$ Recall notation from
\ref{2.2}.

 For $i\in \{j\}_{j=1}^n$ put $r_T(i)$ be
 the number of the row of T $i$ belongs to.

  For $i\in \langle T_1\rangle$
 and $j\in\langle T_2\rangle$ set ${\rm Change}(T,i,j):=(S_1,S_2)$ to be the
 2-column array where elements increase in columns from the top to
 the bottom
 obtained from $T$ by  $\langle S_1\rangle=\{ j\}\cup\langle T_1\rangle\setminus\{i\} $
 and $\langle S_2\rangle=\{i\}\cup\langle T_2\rangle\setminus\{j\}.$ Note that $(S_1,S_2)$ is
 not necessarily a Young tableau.
 Note that in ${\rm Change}(T,i,j)$ we always will write first an element of the first column
 and then an element of the second column.

 Apart from the subsection \ref{5.6}
 in this paper we do not enter the questions
connected to the Robinson-Schensted correspondence. In \cite{MPI} we
had a detailed discussion on the connection between
Robinson-Schensted correspondence and intersections of codimension
1. In connection with this correspondence we need to consider
$s_i:=(i,i+1)$ for $i\ :\ 1\leq i\leq n-1$ which are generators of
$\bS_n$ as the Weyl group of $\G.$

 Translating \cite[\S 3.2, Proposition 3.10]{MPI} to our case we get
 \begin{prop}
 \begin{itemize}
 \item[(i)] Let $T,S\in\bT_\lambda.$ If there exist
 $P\in\bT_\lambda$ and $s_m$ for $m\ :\ 1\leq m\leq n-1$ such that
 ${\rm RS}(T,P)=s_m{\rm RS}(S,P)$ then
 $\codim_{\V_T}\V_T\cap\V_S=1.$
\item[(ii)] Let $T\in \bT_{(n-k,k)}$. If $i\in\langle T_1\rangle$ and
$j=i\pm 1$ is in $\langle T_2\rangle$ are such that $r_T(i)\ne
r_T(j)$ then $S={\rm Change}(T,i,j)$ is a tableau and there exist
$P\in \bT_{(n-k,k)}$ and $s_m$ for $m\ :\ 1\leq m\leq n-1$ such that
${\rm RS}(T,T')=s_m{\rm RS}(S,T').$ In particular
$\codim_{\V_T}(\V_T\cap\V_S)=1.$
\end{itemize}
\end{prop}

\subsection{}\label{5.2ab}
 Given $T\in\bT_{(n-k,k)}$. By Proposition \ref{2.2}
  $\overline\B_{\sigma_T}=\overline\V_T$ where
 $\sigma_T=(i_1,b_1)\ldots(i_k,b_k)$ is defined in \ref{2.2}.

  Note also
\begin{lem} For any $k\leq \frac{n}{2}$ let $T\in\bT_{(n-k,k)}$ and let
$\sigma_T=(i_1,b_1)\ldots(i_k,b_k)$ then for any $(i_s,b_s)\in
\sigma_T$ one has
\begin{itemize}
\item[(i)] if $i_s\ne b_s-1$ then for any $p\:\ i_s<p<b_s$
one has $p\in\{i_q,b_q\}_{q=1}^{s-1};$
\item[(ii)] $b_s-i_s$ is odd.
\end{itemize}
\end{lem}
\begin{proof}
Claim (i) is a straightforward corollary of the definition. Indeed
$i_1=b_1-1$ and for any $s\ :\ 2\leq s\leq k$, we have
$i_s=\max\{m\in\langle T_1\rangle-\{i_j\}_{j=1}^{s-1}\ |\
 m<b_s\}$. So if $i_s\ne b_s-1$ then for any $p\ :\ i_s<p<b_s$ such
 that $p\not\in\langle T_2\rangle$ one has $p=i_q$ for some $q<s.$

Claim (ii) follows from claim (i). Indeed $b_1-i_1=1$. For $s\ :\
2\leq s\leq k$ if $i_s\ne b_s-1$ let $b_p$ be minimal such that
$b_p>i_s.$ For any $q\ :\ p\leq q<s$ one has $b_p>i_s$ so that
$i_p>i_s$ (as maximum) thus
$\{i_q,b_q\}_{q=p}^{s-1}=\{t\}_{t=i_s+1}^{b_s-1}.$ Hence
$b_s-i_s=2(s-p)+1.$
\end{proof}

\subsection{}\label{5.2a}
Translating (i) of the proof of Theorem \ref{thm}  we get the
following. For any $(i_s,b_s)\in\sigma_T$ if
$\{i_q,b_q\}_{q<s}\not\supset\{j\}_{j=1}^{i_s-1}$ let
$m(i_s)=\max\{r\in\{j\}_{j=1}^{i_s-1}\setminus \{i_q,b_q\}_{q<s}\}.$
If either $s=k$ or for any $t>s$ one has $i_t\ne m(i_s)$ then ${\rm
Change}(T,i_s,b_s)$  is a Young tableau and
$\codim_{\V_T}(\V_T\cap\V_{{\rm Change}(T,i_s,b_s)})=1.$

Exactly in the same way translating (ii) of the proof of Theorem
\ref{thm}  we get: For any $(i_s,b_s)\in\sigma_T$ if
$\{i_q,b_q\}_{q>s}\not\supset\{j\}_{j=b_s+1}^n$ let
$m(b_s)=\min\{r\in\{j\}_{j=b_s+1}^{n}\setminus \{i_q,b_q\}_{q>s}\}.$
If either $s=k$ or for any $t>s$ such that $b_t<m(b_s)$ one has
$i_t>i_s$ then
 ${\rm Change}(T,m(b_s),b_s)$ is a Young tableau and
 $\codim_{\V_T}(\V_T\cap\V_{{\rm Change}(T,m(b_s),b_s)})=1.$

Translating (iii) of the proof of Theorem \ref{thm} we get: For
$s<k$ if there exists $t>s$ such that $i_t<i_s$  then ${\rm
Change}(T,i_s,b_s)$ is a Young tableau and
$\codim_{\V_T}(\V_T\cap\V_{{\rm Change}(T,i_s,b_s)})=1.$

Finally (iv) of the proof is translated into: For $s<k$ put
$P((i_s,b_s)):=\{(i_t,b_t)\ :\ i_t>b_s\ {\rm and}\ \forall\ s<q<t:
i_s<i_q,\forall\ p>t\ i_p\not\in\{j\}_{i_s+1}^{i_t-1};\ {\rm and}\
\{i\}_{i=b_s+1}^{i_t-1}\setminus\langle\sigma_T\rangle=\emptyset\}.$
For any $(i_t,b_t)\in P((i_s,b_s))$  ${\rm Change}(T,i_t,b_s)$ is a
Young tableau and $\codim_{\V_T}(\V_T\cap\V_{{\rm
Change}(T,i_t,b_s)})=1.$

The results in the form  provided above are not very clear, however
we can simplify them further noting that (i) and (iii) together
provide us

\begin{prop}\label{prop5.3}
For any $k\leq \frac{n}{2}$ let $T\in\bT_{(n-k,k)}$ and let
$\sigma_T=(i_1,b_1)\ldots(i_k,b_k)$ then for any $s\leq k$ such that
$b_s>2s$ one has ${\rm Change}(T,i_s,b_s)$ is a Young tableau such
that $\codim_{\V_T}(\V_T\cap\V_{{\rm Change}(T,i_s,b_s)})=1.$
\end{prop}

\begin{proof}
Let us show first that the condition $b_s>2s$ is equivalent to the
condition
$$\{j\}_{j=1}^{i_s-1}\setminus\{i_1,b_1,\ldots,i_{s-1},b_{s-1}\}\ne\emptyset.$$
Indeed, $b_s\geq 2s$ always since $a_1,\ldots,a_s$ and
$b_1,\ldots,b_{s-1}$ are smaller than $b_s.$ Now if $b_s=2s$ then
$\{i_t,b_t\}_{t=1}^s=\{j\}_{j=1}^{2s}$ so that
$\{i_1,b_1,\ldots,i_{s-1},b_{s-1}\}=\{j\}_{j=1}^{2s-1}\setminus\{i_s\}$
and in particular $\{j\}_{j=1}^{i_s-1}\subset
\{i_1,b_1,\ldots,i_{s-1},b_{s-1}\} .$ On the other hand if $b_s>2s$
then $\sh(\pi_{1,b_s}(T))=(b_s-s,s)$ where $b_s-s>s$ so that
$i_s=\max\{j\in\{a_j\}_{j=1}^{b_s-s}\setminus\{i_t\}_{t=1}^{s-1}\}$
and there exists $a_j<i_s$ such that $a_j\not\in
\{i_t\}_{t=1}^{s-1}.$ Thus
$a_j\in\{j\}_{j=1}^{i_s-1}\setminus\{i_1,b_1,\ldots,i_{s-1},b_{s-1}\}.$

Let $b_s>2s$ and let
$m(i_s)=\max\{j\in\{j\}_{j=1}^{i_s-1}\setminus\{i_1,b_1,\ldots,i_{s-1},b_{s-1}\}\}.$
One has
\begin{itemize}
\item{} either $s=k$ or for any $t>s$ $i_t>i_s$
 thus the condition of (i) is satisfied;
\item{} or $i_t=m(i_s)$ for some $t>s$, thus the condition of (iii)
is satisfied;
\end{itemize}
In both cases ${\rm Change}(T,i_s,b_s)$ is a Young tableau and
$\codim_{\V_T}(\V_T\cap\V_{{\rm Change}(T,i_s,b_s)})=1.$
\end{proof}
Note that this proposition gives the rule of interchange the
elements of the first and second columns of a tableau when the
element of the first column is smaller than the element of the
second column.

Note also that for $b_s=2a+1$ (i.e. odd) the condition of the
proposition is always satisfied so that for $b_s$ odd ${\rm
Change}(T,i_s,b_s)$ is a Young tableau such that
$$\codim_{\V_T}(\V_T\cap\V_{{\rm Change}(T,i_s,b_s)})=1.$$
\subsection{}\label{5.4} In the same way (ii) and (iv) together
provide
\begin{prop}\label{prop5.4}
For any $k\leq \frac{n}{2}$ let $T\in\bT_{(n-k,k)}$ and let
$\sigma_T=(i_1,b_1)\ldots(i_k,b_k).$ For any $a_s\in T_1$ such that
$a_s-1\in  T_2$ let $B(a_s)=\{b_p\in T_2\ :\ b_p=a_s-1\ {\rm or}\
\{j\}_{b_p+1}^{a_s-1}=\{i_q,b_q\}_{q>p}\}.$ For any $b_p\in B(a_s)$
${\rm Change}(T,a_s,b_p)$ is a Young tableau such that
$$\codim_{\V_T}(\V_T\cap\V_{{\rm Change}(T,a_s,b_p)})=1.$$
\end{prop}
\begin{proof}
Note that for any $a_s,b_p$ such that $a_s>b_p$ one has $S={\rm
Change}(T,a_s,b_p)$ is a Young tableau. Let us show that for $b_p\in
B(a_s)$
 $(b_p,a_s)\in \sigma_S$ so that $T={\rm Change}(S,b_p,a_s)$
 from proposition \ref{prop5.3}.
 Then the result follows. Note that if $a_s-1\in T_1$ then
 $(a_s-1,a_s)\in\sigma_S$ so that $\V_T\cap\V_S$ cannot be of
 codimension 1 by proposition \ref{prop5.3}. Thus, the condition
 $a_s-1\in T_2$ is necessary.

By definition of $\sigma_S$ this is obvious for $b_p=a_s-1.$ Now
assume
$$\{j\}_{b_p+1}^{a_s-1}=\{i_q,b_q\}_{q>p}.\eqno{(*)}$$
Let us denote $b_t=a_s-1.$ Then
$S_2=(b_1,\ldots,b_{p-1},b_{p+1},\ldots,b_t,a_s,b_{t+1},\ldots).$
Thus
 $(i_q,b_q)\in \sigma_S$ iff $(i_q,b_q)\in\sigma_T$ for any $q\ :\
 q\leq p-1.$ Further note that for any $q\ :\ p+1<q\leq t$ by the
 conditions $(*)$ $i_q>b_p$ so that again $(i_q,b_q)\in \sigma_S$ iff
 $(i_q,b_q)\in\sigma_T$
 for any $q\ :\ p+1<q<t.$ Moreover by $(*)$ $b_p=\max\{j\in \langle
 S _1\rangle\setminus\{i_q\}_{q=1,q\ne p}^t\ :\ j<a_s\}$ thus
 $(b_p,a_s)\in \sigma_S.$
Consider $\sigma_S$ and note that $a_s>2p$ (since $a_s>b_p\geq 2p$)
so that the conditions of proposition \ref{prop5.3} are satisfied.
Thus $T={\rm Change}(S,b_p,a_s)$ is such that
$\codim_{\V_T}(\V_T\cap\V_S)=1.$
\end{proof}
Note that this proposition gives the rule of interchange the
elements of the first and second columns of a tableau when the
element of the first column is greater than the element of the
second column.
\subsection{}\label{5.5}

Note also that propositions \ref{prop5.3}, \ref{prop5.4} provide us
\begin{cor}
If $\codim_{\V_T}(\V_T\cap\V_{{\rm Change}(T,i,j)})=1$ then $j-i$ is
odd.
\end{cor}
\begin{proof}
If ${\rm Change}(T,i,j)$ is from proposition \ref{prop5.3} then
$(i,j)\in \sigma_T$.  If ${\rm Change}(T,i,j)$ is from proposition
\ref{prop5.4} then $(j,i)\in\sigma_{{\rm Change}(T,i,j)}.$ In both
cases $j-i$ is odd by lemma \ref{5.2ab}.
\end{proof}
 \vskip 0.2cm
\subsection{}\label{5.6} In \cite{MPI} we promised to show in this paper
 \begin{prop}\label{prop5.6} Let $T,S\in\bT_{(n-2,2)}$ then
 $\codim_{\V_T} \V_T\cap\V_S=1$ if and only if there exists
$P\in\bT_{(n-2,2)}$ such that ${\rm RS}(T,P)=s_i {\rm RS}(S,P).$
\end{prop}
\begin{proof}

We have only two possibilities for $T$. Either $T_2=\{a,b\}$ where
$b-a\geq 2,$ or $T_2=\{a,a+1\}.$

Note that it is enough to check only ${\rm Change}(T,i,j)$ where
$i<j$ that is they are from Proposition \ref{prop5.3} since if $i>j$
then for $T'={\rm Change}(T,i,j)$ one has $T={\rm Change}(T',j,i)$
where $j<i.$

(a) If $T_2=\{a,b\}$ where $b-a\geq 1$ then
$\sigma_T=(a-1,a)(b-1,b).$ By Proposition \ref{prop5.3} one has two
possible tableaux of form ${\rm Change}(T,i,j)$ where $i<j.$
\begin{itemize}
\item[(i)] If $a\geq 3$ one has the intersection of codimension 1 with
$\V_S$ where
$$S={\rm Change}(T,a-1,a)=\begin{array}{lc}
1&a-1\\
\vdots&b\\
\vdots&\\
\end{array}$$

\item[(ii)] If $b\geq 5$ one has the intersection of codimension 1 with
$\V_S'$ where
$$S={\rm Change}(T,b-1,b)=\begin{array}{lc}
1&a\\
\vdots&b-1\\
\vdots&\\
\end{array}$$
\end{itemize}
and by Proposition \ref{5.2}(ii) in both cases there exist
$P\in\bT_{(n-2,2)}$ and $s_m$ for $m\ :\ 1\leq m\leq n-1$ such that
${\rm RS}(T,P)=s_m {\rm RS}(S,P)$ (resp. ${\rm RS}(T,P)=s_m {\rm
RS}(S',P)$

 (b) If $T_2=\{a,a+1\}$ then
$\sigma_T=(a-1,a)(a-2,a+1).$ By Proposition \ref{prop5.3} one has
two possible tableaux of form ${\rm Change}(T,i,j)$ where $i<j.$
\begin{itemize}
\item[(i)] Since $a\geq 3$ one always has the intersection of codimension 1 with
$\V_S$ where
$$S={\rm Change}(T,a-1,a)=\begin{array}{ll}
1&a-1\\
\vdots&a+1\\
\vdots&\\
\end{array}$$
 and by Proposition \ref{5.2}(ii) there exist $P\in\bT_{(n-2,2)}$ and $s_m$
  for $m\ :\ 1\leq m\leq n-1$ such
that ${\rm RS}(T,P)=s_m {\rm RS}(S,P).$

\item[(ii)] If $b+1>4$ then one has the intersection of codimension 1 with
$\V_S'$ where
$$S'={\rm Change}(T,a-2,a+1)=\begin{array}{lc}
1&a-2\\
\vdots&a\\
\vdots&\\
\end{array}$$
Let us construct $P$ and $s_m$ such that ${\rm RS}(T,P)=s_m {\rm
RS}(S',P).$ Indeed, let us take
$$P=\begin{array}{cc}
1&n-a+2\\
\vdots&n-a+3\\
\vdots&\\
n-a&\\
n-a+3&\\
\vdots&\\
n&\\
\end{array}$$
One can easily check that
 $${\rm RS}(T,P)=[n,\ldots,a+2,a-1,{\bf
a-2},a+1,a,{\bf a-3},a-4,\ldots,1]$$
 and correspondingly
 $${\rm
RS}(S',P)=s_{a-3}{\rm RS}(T,P)=[n,\ldots,a+2,a-1,{\bf
a-3},a+1,a,{\bf a-2},a-4,\ldots,1]$$
\end{itemize}
\end{proof}

\subsection{}\label{5.7}
We would like to connect this section to the results of \cite{MPI}.
We start with a few notes.
\begin{rem}
\begin{itemize}
\item[(i)]  Note that the last case of the proof of Proposition \ref{prop5.6}
is the only case among 4 cases we have considered when $S'$
 is not Vogan's
 $\T_{\alpha,\beta}(T)$ (for the definition and details of
 $\T_{\alpha,\beta}$ cf. \cite[\S 3]{MPI}).

\item[(ii)] As we have noted in \cite[4.2]{MPI} for $\bT_{(n-1,1)}$ one
has $\V_T\cap\V_S$ is of codimension 1 iff $S=\T_{\alpha,\beta}(T)$
so that in particular $S={\rm Change}(T,i,j)$ where $j=i\pm 1$. Thus
by Proposition \ref{5.2} there exist $P\in\bT_{(n-1,1)}$ and $s_m$
such that ${\rm RS}(T,P)=s_m{\rm RS}(S,P).$ As well we have shown in
\cite[5.8]{MPI} that for $\bT_{(n-k,k)}$ where $k\geq 3$ there exist
$S,T\in\bT_{(n-k,k)}$ such that  $\codim_{\V_T}(\V_S\cap\V_T)=1$,
however there is no $P\in \bT_{(n-k,k)}$ such that ${\rm
RS}(T,P)=s_m{\rm RS}(S,P)$ for some $m\ :\ 1\leq m\leq n-1.$
\end{itemize}
 \end{rem}

 Let us finish this article by a speculation about the
intersection components in codimension one in the Springer fiber.
The hook case and the two column case are two extreme cases in the
following sense: For all the nilpotent orbits of the given rank $k$
the orbit $\lambda=(k,1,\ldots)$ is the most non-degenerate and the
orbit $\lambda=(2,\cdots,2,1,\ldots)$ (with dual partition
$(n-k,k)$) is the most degenerate, in the following sense
$\overline\Or_{(k,1,\ldots)}\supset
\overline\Or_{\mu}\supset\Or_{(2,\ldots,2,1,\cdots)}$.

The intersections of the components of the Springer fiber in hook
case was studied by J.A. Vargas \cite{Var}. In particular he showed
that all the intersections are irreducible. However this is not true
in general, in \cite{MPI} we showed that the intersections are
reducible in general and have components of different dimensions.
Nevertheless by Theorem \ref{thm} the intersection in codimension
one of two irreducible components is always irreducible in the two
column case. So the following conjecture seems to be plausible:

\begin{conj}
The intersection in codimension one of two irreducible components of
a Springer fiber of type $A$ is irreducible.
\end{conj}

\vskip 0.2 cm {\bf Acknowledgements.}  The second author would like
to express his gratitude to L\^{e} D\~{u}ng Tr\'{a}ng, B. Dubrovin
and C.S. Seshadri for the invitation to the Abdus Salam
International Centre for Theoretical Physics (I.C.T.P), at the
Department of Mathematics International School for Advanced Studies
(S.I.S.S.A / I.S.A.S) in Trieste, and at the Chennai Mathematical
Institute (C.M.I) in India, where this work was done. He would like
also thank all these institutes for their kind hospitality and
support.

\bigskip

\centerline{ INDEX OF NOTATION}

\noindent Symbols and notions appearing frequently are given below
in order of appearance.

\ref{1.1} $\bG,\ \G,\ g.u, \nil, \bB,\ \bS_n,\ \Or_x,\ \V_T;$

\ref{1.2} $\sh(x),\ \sh(T),\ \bT_{(n-k,k)};$

\ref{1.3} $\X^2,\ \B_x,\ \Bb^2,\ \Bb_{(n-k,k)},\ \bS_n^2,\
\bS_n^2(k);$

\ref{2.1} $N_\sigma, \B_\sigma,\ L(\sigma);$

\ref{2.2} $T_1,\ T_2, \langle T_i\rangle, \sigma_T;$

\ref{2.3} canonical form of $\sigma$,
$q_s(\sigma)=q_{(i_s,j_s)}(\sigma);$

\ref{2.4} $\pi_{i,j},\ R_u,\ R_\sigma,\ \bR_n^2;$

\ref{2.5} $\preceq;$

\ref{2.6} $R_{\sigma,\sigma'};$

\ref{4.4} $\sigma_{(i,j)}^-;$

\ref{4.5} $D(\sigma),\ A(\sigma);$

\ref{4.6a} $\langle \sigma\rangle,\ \langle \sigma\rangle^c;$

\ref{4.6b} $F_{(p\rightarrow q)}(\sigma),\ F_{(p\rightleftarrows
q)}(\sigma),$ minimal elementary transformation;

\ref{4.61} $\sigma_{i_s\dar},\ \sigma_{\uar i_s};$

\ref{4.62} $\sigma_{j_s\rar},\ \sigma_{\lar j_s};$

\ref{4.63} $C_{i_t\looparrowright}(\sigma),\ C_{\looparrowleft
i_t}(\sigma);$

\ref{4.64} $C_{i_s\uar\dar}(\sigma),\ C_{\dar\uar i_s}(\sigma);$

\ref{5.2} $r_T(i),\ {\rm Change}(T,i,j);$


\begin{thebibliography}{mmm}

\bibitem[Hot]{Hot}
R. Hotta, On Joseph's construction of Weyl group representations,
Tohoku Math. J. {\bf 36} (1984), 49-74.

\bibitem[K-L]{KL}
D. Kazhdan and G. Lusztig,  Representations of Coxeter groups and
Hecke algebras, Invent. Math.  {\bf 53} (1979) 165-184.

\bibitem[Mel1]{Msmith} A.Melnikov,  Orbital varieties in $\Sl_n$  and the Smith conjecture,
                 J. of Algebra, {\bf 200} (1998),  1-31.

\bibitem[Mel2]{Mx2}    A.Melnikov, B-orbits in solutions to the equation $X^2=2$ in triangular
matrices,   J. of Algebra, {\bf 223} (2000),  101-108.

\bibitem[Mel3]{Mel1}
A. Melnikov, Description of B-orbit closures of order 2 in
upper-triangular matrices,   Transform. Groups  {\bf 11}  (2006),
no. 2, 217--247.

\bibitem[MP1]{M-P} A. Melnikov, N.G.J. Pagnon,
On intersections of orbital varieties and components of Springer
fiber, J. Algebra, {\bf 298} (2006),  1-14.

\bibitem[MP2]{MPI} A. Melnikov, N.G.J. Pagnon, Intersections of
components of a Springer fiber for the hook and two column cases,
Arxiv, math.RT/0607673.

\bibitem[Spa1]{Spa0}
N. Spaltenstein,  The fixed point set of a unipotent transformation
on the flag manifold, Indag. Math. {\bf 38} (1976), 452-456.

\bibitem[Spa2]{Spa}
N. Spaltenstein,  On the fixed point set of a unipotent element on
the variety of Borel subgroups, Topology {\bf 16} (1977), 203-204.

\bibitem[Ste]{Ste2}
R. Steinberg,  On the Desingularisation of the Unipotent Variety,
Invent. Math.  {\bf 36}  (1976), 209-224.

\bibitem[Var]{Var}
J.A. Vargas,  Fixed points under the action of unipotent elements of
${\rm {SL}_{n}}$ in the flag variety, Bol. Soc. Mat. Mexicana 2{\bf
4} (1979), 1-14.
\end{thebibliography}
\end{document}